\documentclass{article}
\usepackage{amssymb,amscd,amsfonts,amsbsy}
\def\ring{\mathaccent"0017}
\def\ring{\mathaccent"0017}

\newtheorem{proposition}{Proposition}[section]
\newtheorem{theorem}[proposition]{Theorem}
\newtheorem{lemma}[proposition]{Lemma}
\newtheorem{corollary}[proposition]{Corollary}
\newtheorem{remark}{Remark}[section]

\title{{Analytic criteria in the qualitative spectral analysis of the Schr\"odinger operator}
\thanks{
Supported in part by  NSF Grant DMS-0500029
}}

\author{Vladimir Maz'ya
\\*[10pt]
\emph{\small Ohio State University, Columbus, OH, USA}
\\
\emph{\small  University of Liverpool, Liverpool, UK}
\\
\emph{\small  Link\"oping University, Sweden}
\\*[10pt]
\emph{Dedicated to Barry Simon on the occasion of  his 60th birthday}}
\date{}

\begin{document}
\maketitle
\vspace{-.5cm}
\begingroup

\narrower\noindent
{\it Abstract.}  A number of topics in the qualitative spectral analysis of the  Schr\"odinger operator $-\Delta + V$ are surveyed. In particular, some old and new results concerning the positivity and semiboundedness of this operator  as well as  the structure of different parts of its spectrum are considered. The attention is focused on conditions both necessary and sufficient, as well as on their sharp corollaries.
\endgroup

{\renewcommand{\baselinestretch}{.1} \scriptsize
\tableofcontents
}
\newpage

%\medskip\noindent
%{ {\small\it Key words}: }{ Schr\"odinger operator, qualitative spectral analysis, Wiener's capacity, pointwise multipliers }

\section{Introduction}

 The purpose of the present  article is to survey various analytic results concerning the Schr\"odinger operator $-\Delta + V(x)$  
 obtained during the last half-century, including several quite recently.  
  In particular, the positivity and semiboundedness of this operator  as well as  the structure of different parts of its spectrum are among the topics touched upon.  
  
In the choice of material I aim mostly at results of final character, i.e. at 
simultaneously necessary and sufficient conditions,  and their consequences, the best possible in a sense.   
Another  motivation for the inclusion of any particular  problem  in this survey  is my own 
involvement in its solution.

Naturally, the selection of topics is far from 
exhaustive. For instance, conditions for the essential self-adjointness of the Schr\"odinger
operator and for the absence of eigenvalues at positive energies, as well as properties of eigenfunctions and bounds on the number of eigenvalues,  are not considered here. These and many other  themes of qualitative spectral analysis have been discussed in the comprehensive surveys  \cite {[Sim4]}, \cite {[BMS]}, \cite {[Dav1]}, and, of course, in the classical Reed and Simon's treatise  \cite{[RS1]}-\cite{[RS4]}.  However, the information collected in what follows has been  quoted rarely if ever  in the quantum mechanics literature.

 Due to space limitation,  I write only about the linear Schr\"odinger operator with
electric potential, although some of the results below can be
modified and applied to  nonlinear, magnetic,
 relativistic Schr\"odinger operators  and even general elliptic operators with variable coefficients.
 For  the same reason, proofs are supplied just in a few cases, mostly
 when a source does not seem readily available.

\section{Wiener capacity}

The capacity of a set in $\Bbb{R}^n$ will be  met frequently in the present article. This notion appeared first in electrostatics and was introduced to mathematics by N. Wiener in the 1920s. Since then several generalizations and modifications of  Wiener's capacity appeared: Riesz, Bessel, polyharmonic capacities, $p$-capacity and others.  They are of use in potential theory, probability, function theory and partial differential equations. The capacities provide adequate terms to describe sets of discontinuities of Sobolev functions,  removable singularities of solutions to partial differential equations, sets of uniqueness for analytic functions, regular boundary points in the Wiener sense, divergence sets for trigonometric series, etc. (see, for example, \cite{[AH]}, \cite {[Car]}, \cite{[Lan]}, \cite {[MH]}). Some applications of Wiener's and other capacities to the theory of the Schr\"odinger operator are presented in this survey.

Let $\Omega$ be an open set in $\Bbb{R}^n$ and let $F$ be a compact subset of $\Omega$. The Wiener (harmonic) capacity of $F$ with respect to $\Omega$ is defined as the number

\begin{equation}\label{1.18}
{\rm cap}_\Omega F = \inf \Bigl\{ \int_\Omega| \nabla \, u|^2 dx: \,\, u\in C^\infty_0(\Omega), \,\, u\geq 1\,\, {\rm on}\,\, F\Bigr\}.
\end{equation}
 We shall use the simplified notation ${\rm cap} F$ if
$\Omega = \Bbb{R}^n$. The capacity ${\rm cap}_\Omega F$ can be 
  defined equivalently 
as the least upper bound of $\nu(F)$ over the set of all measures
$\nu$ supported by $F$ and satisfying the condition
$$\int_\Omega G(x, y)\, d\nu(y)\leq 1,$$
where $G$ is the Green function of the domain $\Omega$. If $\Omega = \Bbb{R}^3$ then it is just the electrostatic capacity of $F$.

It follows from the definition (\ref{1.18}) that the capacity is a nondecreasing function of $F$ and a nonincreasing one  of $\Omega$.  We have Choquet's inequality
$${\rm cap}_\Omega(F_1\cap F_2) + {\rm cap}_\Omega(F_1\cup F_2) \leq {\rm cap}_\Omega \, F_1 + {\rm cap}_\Omega \, F_2 $$
for arbitrary compact sets $F_1$ and $F_2$ in $\Omega$ \cite{[Cho]}. It is easy to check that the Wiener capacity is continuous from the right. This means that for each $\varepsilon >0$ there exists a neighborhood $G$, $F\subset G \subset\overline{G} \subset \Omega$ such that for each compact set $F_1$ with $F\subset F_1 \subset G$ the inequality
$${\rm cap}_\Omega\, F_1\leq {\rm cap}_\Omega \, F +\varepsilon$$
holds.

Let $E$  be an arbitrary subset of $\Omega$. The inner and the outer capacities are defined as numbers
$$\underline{\rm cap}_\Omega\, E =\mathop{\hbox
{\rm sup}}_{F\subset E} \,{\rm cap}_\Omega \, F,\qquad F\,\, {\rm compact}\,\, {\rm in}\,\,\Omega,$$

$$\overline{\rm cap}_\Omega\, E =\mathop{\hbox
{\rm inf}}_{G\supset E} \,{\rm cap}_\Omega\, G,\qquad G\,\, {\rm open}\,\, {\rm in}\,\,\Omega.$$

It follows from the general Choquet theory that for each Borel set both capacities coincide \cite{[Cho]}. Their common value is called the Wiener (harmonic) capacity and will be denoted by  ${\rm cap}_\Omega\, E$.

By $v_n$ we denote the volume of the unit ball in $\Bbb{R}^n$ and 
let ${\rm mes}_n F$ stand for the $n$-dimensional Lebesgue measure of $F$.
By the classical isoperimetric inequality, the following isocapacitary inequalities hold (see \cite{[Maz7]}, Sect. 2.2.3)
\begin{equation}\label{a1}
{\rm cap}_\Omega \, F \geq nv_n^{2/n}(n-2)\, |({\rm mes}_n\, \Omega)^{(2-n)/n} - ({\rm mes}_n\, F)^{(2-n)/n}|^{-1}\quad {\rm if}\,\, n>2
\end{equation}
 and
\begin{equation}\label{a2}
{\rm cap}_\Omega F \geq 4\pi\Bigl(\log\frac{{\rm mes}_2\, \Omega}{{\rm mes}_2\, F}\Bigr)^{-1}\quad {\rm if}\,\, n=2.
\end{equation}
In particular, if $n>2$ then
\begin{equation}\label{a3}
{\rm cap} F \geq nv_n^{2/n}(n-2)({\rm mes}_n\, F)^{(n-2)/n}.
\end{equation}
If $\Omega$ and $F$ are concentric balls, then the three preceding estimates come as identities.

Using Wiener's capacity, one can obtain two-sided estimates for the best constant in the Friedrichs inequality
\begin{equation}\label{a4}
\|u\|_{L_2(B_1)} \leq	 C\, \|\nabla u\|_{L_2(B_1)},
\end{equation}
where $B_1$ is a unit open ball and $u$ is an arbitrary function in $C^\infty(\overline{B_1})$ vanishing on a compact subset $F$ of $\overline{B_1}$. 

\begin{proposition}\label{P1}
{\rm \cite{[Maz2]}}  The best constant $C$ in {\rm (\ref{a4})} satisfies
\begin{equation}\label{a5}
C\leq c(n)\bigl( {\rm cap}\, F\bigr)^{-1/2},
\end{equation}
where $c(n)$ depends only on $n$. 
\end{proposition}

It is shown in \cite{[MShu2]} that (\ref{a5}) holds with $c(n) = 4v_n n^{-1}(n^2-2)$. 
Proposition \ref{P1} has the following partial converse. 

\begin{proposition}\label{P2}
  Let
\begin{equation}\label{a6}
{\rm cap}\, F \leq \gamma\, {\rm cap} \, B_1,
\end{equation}
where $\gamma\in (0,1)$. Then any  constant $C$  in {\rm (\ref{a4})} satisfies
\begin{equation}\label{a7}
C\geq c(n,\gamma)\bigl( {\rm cap}\, F\bigr)^{-1/2}.
\end{equation}
\end{proposition}

This assertion was proved  in \cite{[Maz2]} (see also \cite{[Maz7]}, Ch.\,10) with  a sufficiently small $\gamma =\gamma(n)$. The present stronger version is essentially contained  in \cite{[MShu1]}. Both propositions   were used in \cite{[KMS]} to derive a discreteness  of spectrum  criterion for the magnetic Schr\"odinger operator.

\section{Equality of the minimal and maximal Dirichlet Schr\"odinger forms}

Let $\Bbb{V}$ be a nonnegative Radon measure in $\Omega$. Consider the quadratic form
\begin{equation}\label{b1}
Q[u,u] = \int_\Omega |\nabla\, u|^2\, dx + \int_\Omega |u|^2\, \Bbb{V}(dx).
\end{equation}
The closure of $Q$ defined on the set
$$\{ u\in C^\infty(\Bbb{R}^n) :\, {\rm supp}\, u\subset \Omega, \,\, Q[u,u] < \infty\},$$
which may contain functions with non-compact support, 
will be denoted by $Q_{\rm max}$. Another quadratic form  $Q_{\rm min}$ is introduced as the closure of $Q$ defined on $C^\infty_0(\Omega)$.

The question as to when the equality $Q_{\rm max}=Q_{\rm min}$ holds has been raised by T.\,Kato \cite{[Ka3]}.  B.\,Simon studied a similar question concerning the magnetic Schr\"odinger operator in $\Bbb{R}^n$ \cite{[Sim2]}. The following necessary and sufficient condition for this equality was obtained in \cite{[CMaz]}.

\begin{theorem}\label{Th1}
{\rm (i)} If either $n\leq 2$ or $n>2$ and ${\rm cap}(\Bbb{R}^n\backslash\Omega) = \infty$, then $Q_{\rm max}=Q_{\rm min}$.

{\rm (ii)} Suppose that $n>2$ and ${\rm cap}(\Bbb{R}^n\backslash\Omega) < \infty$. Then $Q_{\rm max}=Q_{\rm min}$ if and only if
\begin{equation}\label{1}
\Bbb{V}(\Omega\backslash F) = \infty \quad {\rm for} \,\, {\rm every}\,\,  {\rm closed}\,\, {\rm set}\,\, F\subset\Omega\,\, {\rm with}\,\,  {\rm cap}F <\infty.
\end{equation}
\end{theorem}

\noindent
{\bf Example.}  Consider the domain $\Omega$ complementing the 
infinite funnel 
$$\{ x=(x', x_n): \, x_n \geq 0, \, |x'|\leq f(x_n)\},$$
where $f$ is a continuous and  decreasing function on $[0,\infty)$ subject to   $f(t) \leq c f(2t)$. One can show that ${\rm cap}\, (\Bbb{R}^n\backslash \Omega) <\infty$ if and only if the function $f(t)^{n-2}$ for $n>3$ and the function $(\log(t/f(t)) )^{-1}$ for $n=3$ are integrable on $(1,\infty)$.

\begin{remark}\label{rem1}
{\rm
 Note that the equality $Q_{\rm max}=Q_{\rm min}$ is equivalent to (\ref{1}) in the particular case $\Omega = \Bbb{R}^n$. Assume that $\Bbb{V}(dx) = V(x)\, dx$, where $V$ is a positive function in $L_{\rm loc}^1(\Bbb{R}^n)$ locally bounded away from zero. According to \cite{[CMaz]}, condition (\ref{1}) is necessary for the essential self-adjointness of the operator $V^{-1}\Delta$ on $L^2(V\, dx)$ with domain $C^\infty_0(\Bbb{R}^n)$.
}
\end{remark}

\section{Closability of quadratic forms}

Let $\Omega$ be an open set in $\Bbb{R}^n$ and let $\Bbb{V}$ be a nonnegative Radon measure in $\Omega$. By $\ring L^1_2(\Omega)$ we denote the completion of $C^\infty_0(\Omega)$ in the norm $\|\nabla\, u\|_{L_2(\Omega)}$.

One says that the quadratic form
\begin{equation}\label{N1}
\|u\|^2_{L_2(\Omega,\Bbb{V}(dx))} : = \int_\Omega |u|^2\, \Bbb{V}(dx)
\end{equation}
defined on $C^\infty_0(\Omega)$ is closable in $\ring L^1_2(\Omega)$ if any Cauchy sequence in the space $L_2(\Omega,\Bbb{V}(dx))$ converging to zero in $\ring L^1_2(\Omega)$
has zero limit in $L_2(\Omega,\Bbb{V}(dx))$. We call the measure  $\Bbb{V}$ absolutely continuous with respect to the harmonic capacity if the equality ${\rm cap}\, E =0$, where $E$ is a Borel subset of $\Omega$, implies $\Bbb{V}(E) =0$. $\qquad\square$

Let, for example, $\Bbb{V}$ be the $\varphi$-Hausdorff measure in $\Bbb{R}^n$, i.e.
$$\Bbb{V}(E) = \mathop{\hbox
{\rm lim}}_{\varepsilon\to +0} \mathop{\hbox
{\rm inf}}_{\{{\cal B}^{(j)}\}} \sum_i\varphi(r_i),$$
where $\varphi$ is a nondecreasing positive continuous function on $(0,\infty)$ and $\{{\cal B}^{(j)}\}$ is any covering of the set $E$ by open balls ${\cal B}^{(i)}$ with radii $r_i<\varepsilon$. It is well known  that this measure is absolutely continuous with respect to the Wiener capacity if
$$\int_0^\infty \varphi (t)\, t^{1-n}\, dt <\infty,$$
and that  this condition is sharp in a sense (see \cite{[Car]}).

Needless to say, any measure $\Bbb{V}$ absolutely continuous with respect to  the $n$-dimensional Lebesgue measure is absolutely continuous with respect to the capacity. $\qquad\square$

The following result has been obtained in \cite{[Maz3]}.

\begin{theorem}\label{Th2}
The quadratic form (\ref{N1}) is closable in $\ring L^1_2(\Omega)$  if and only if $\Bbb{V}$ is absolutely continuous with respect to the Wiener capacity.
\end{theorem}

The closability of (\ref{N1}) in $\ring L^1_2(\Omega)$ is necessary for the Schr\"odinger operator $-\Delta - \Bbb{V}$, formally associated with the form
\begin{equation}\label{b0}
S[u,u]:= \int_\Omega |\nabla\, u|^2\, dx - \int_\Omega |u|^2\, \Bbb{V}(dx),\quad u\in C_0^\infty(\Omega),
\end{equation}
to be well-defined. $\qquad\square$

Dealing with the Schr\"odinger operator formally given by the expression $-\Delta + \Bbb{V}$ and acting in $L_2(\Omega)$ we need the following notion of a form closable in $L_2(\Omega)$.

By definition, the quadratic form $Q$ defined on $C^\infty_0(\Omega)$ by (\ref{b1})
 is closable in $L_2(\Omega)$ if any Cauchy sequence in the norm $Q[u,u]^{1/2}$ which tends to zero in $L_2(\Omega)$ has zero limit in the norm $Q[u,u]^{1/2}$. This notion is equivalent to the lower semicontinuity of $Q$ in $L_2(\Omega)$ \cite{[Sim1]}.

\begin{theorem}\label{Th3}
{\rm (see\cite{[Maz7]}, Sect. 12.4, 12.5)} The form $Q$ defined on $C^\infty_0(\Omega)$ is closable in $L_2(\Omega)$ if and only if $\Bbb{V}$ is absolutely continuous with respect to the Wiener capacity.
\end{theorem}

We  assume in the sequel that the measure $\Bbb{V}$ is absolutely continuous with respect to the Wiener capacity.

\section{Positivity of the Schr\"odinger operator with negative potential}

The next result was obtained in  \cite{[Maz1]}, \cite{[Maz3]} (see also \cite{[Maz7]}, Th. 2.5.2).

\begin{theorem}\label{Th4}
Let $\Omega$  be an open set in $\Bbb{R}^n$, $n\geq 1$,  and let $\Bbb{V}$ be a nonnegative Radon measure in $\Omega$.
The inequality
\begin{equation}\label{c1}
\int_\Omega|u|^2\, \Bbb{V}(dx) \leq \int_\Omega|\nabla\, u|^2\, dx
\end{equation}
holds for every $u\in C^\infty_0(\Omega)$ provided
\begin{equation}\label{c2}
\frac{\Bbb{V}(F)}{{\rm cap}_\Omega F} \leq \frac{1}{4}
\end{equation}
for all compact sets $F\subset \Omega$.

A necessary condition for (\ref{c1}) is
\begin{equation}\label{c3}
\frac{\Bbb{V}(F)}{{\rm cap}_\Omega F} \leq 1,
\end{equation}
where $F$ is an arbitrary compact subset of $\Omega$.
\end{theorem}

\begin{remark}\label{2a}
{\rm
In inequalities (\ref{c2}), (\ref{c3}) and elsewhere in similar cases, we tacitly assume that vanishing of denominator implies vanishing of numerator, and we may choose any appropriate value of the ratio.
}
\end{remark}

\begin{remark}\label{2}
{\rm
 The necessity of (\ref{c3}) is trivial. The proof of sufficiency of (\ref{c2})
  in \cite{[Maz3]}, \cite{[Maz7]} shows that  (\ref{c2}) implies (\ref{c1}) even without the requirement $\Bbb{V}\geq 0$, i.e. for an arbitrary locally finite real valued charge in $\Omega$. 
}
\end{remark}

\begin{remark}\label{2b}
{\rm The formulation and the proof of Theorem \ref{Th4} do not change if we assume that $\Omega$ is
 an open subset of an arbitrary  Riemannian manifold.
}
\end{remark}

Theorem \ref{Th4} immediately gives the following criterion.

\begin{corollary}\label{Cor1}
The trace inequality
\begin{equation}\label{c5}
\int_\Omega|u|^2\, \Bbb{V}(dx) \leq C \int_\Omega|\nabla\, u|^2\, dx
\end{equation}
holds for every $u\in C^\infty_0(\Omega)$ if and only if
$$\sup_{F\subset\Omega}\frac{\Bbb{V}(F)}{{\rm cap}_\Omega F} <\infty.$$
\end{corollary}

The bounds $1/4$ and $1$ in (\ref{c2}) and (\ref{c3}) are sharp. The gap between these sufficient and necessary conditions is the same as in Hille's non-oscillation criteria for the operator
$$-u'' -\Bbb{V}\, u, \qquad \Bbb{V}\geq 0,$$
on the positive semiaxis $\Bbb{R}^1_+$:
\begin{equation}\label{c4}
x\, \Bbb{V}((x,\infty)) \leq 1/4 \qquad {\rm and}\quad x\, \Bbb{V}((x,\infty)) \leq 1
\end{equation}
for all $x\geq 0$ \cite{[Hil]}. By the way,  conditions (\ref{c4}) are  particular cases of (\ref{c2}) and (\ref{c3}) with $n=1$ and $\Omega = \Bbb{R}^1_+$. $\qquad\square$

Combining Theorem \ref{Th4} with isocapacitary inequalities (\ref{a1})--(\ref{a3}), we arrive at sufficient conditions for (\ref{c1}) whose formulations involve no capacity. For example, in the two-dimensional case,  (\ref{c1}) is guaranteed by the inequality
$$ \Bbb{V}(F) \leq \frac{4\pi}{\log\displaystyle{\frac{{\rm mes}_2\, \Omega}{{\rm mes}_2\, F}}}.$$
The sharpness of this condition can be easily checked by analyzing the well known Hardy-type inequality
$$\int_\Omega\frac{|u(x)|^2}{|x|^2(\log|x|)^2}\, dx \leq 4\int_\Omega |\nabla\, u(x)|^2\, dx,$$
where $u\in C^\infty_0(\Omega)$ and $\Omega$ is the unit disc. $\qquad\square$

The sufficiency of (\ref{c2}) for the inequality (\ref{c1}) can be directly obtained  from the following more precise result.

\begin{theorem}\label{Th5}
{\rm \cite{[Maz8]}} Let $n\geq 1$ and let $\nu$  be a
non-decreasing function on $(0,\infty)$ such that
$$\Bbb{V}(F) \leq \nu({\rm cap}_\Omega F)$$
for all compact subsets $F$ of $\Omega$. If
$$\int_0^\infty | v(\tau)|^2\, |d\nu(\tau^{-1})| \leq \int_0^\infty |v'(\tau)|^2 d\tau$$
for all absolutely continuous $v$ with 
$v' \in L_2(0,\infty)$ and $v(0) = 0$, then {\rm (\ref{c1})} holds.
\end{theorem}

One example of the application of this assertion is the following improvement of
 the Hardy inequality which cannot be deduced from (\ref{c2}):

 $$\int_{\Bbb{R}^n\backslash B_1} u^2\Bigl( 1+ \frac{1}{(n-2)^2(\log|x|)^2}\Bigr)\frac{dx}{|x|^2} \leq \frac{4}{(n-2)^2}\int_{\Bbb{R}^n\backslash B_1} |\nabla u|^2\, dx,$$

 \noindent
 where $u\in C^\infty_0(\Bbb{R}^n\backslash\overline{B}_1)$ and $n>2$. $\qquad\square$

 Here is a dual assertion to Theorem \ref{Th5} which is stated in terms of the Green function $G$ of $\Omega$ and does not depend on the notion of capacity.

\begin{theorem}\label{Th6a}
Let $\Bbb{V}_F$ be the restriction of the measure $\Bbb{V}$ to a compact set $F\subset\Omega$.
Inequality {\rm (\ref{c1})} holds for every $u\in C^\infty_0(\Omega)$ provided
\begin{equation}\label{d1}
\int_\Omega\!\int_\Omega G(x,y)\, \Bbb{V}_F(dx)\, \Bbb{V}_F(dy) \leq \frac{1}{4} \Bbb{V}(F)
\end{equation}

\noindent
for all  $F$. Conversely, inequality {\rm (\ref{c1})} implies
\begin{equation}\label{d2}
\int_\Omega\!\int_\Omega G(x,y)\, \Bbb{V}_F(dx)\, \Bbb{V}_F(dy) \leq \,  \Bbb{V}(F).
\end{equation}
\end{theorem}

{\bf Sketch of the proof.} Let $u$ be a nonnegative function in $C^\infty_0(\Omega)$ such that $u\geq 1$ on $F$. Then
$$\Bbb{V}(F) \leq \int_\Omega u(x)\, \Bbb{V}_F(dx) \leq \Bigl(\int_\Omega\!\int_\Omega G(x,y)\, \Bbb{V}_F(dx)\, \Bbb{V}_F(dy)\Bigr)^{1/2}\|\nabla\, u\|_{L_2(\Omega)}$$
which in combination with (\ref{d1}) gives (\ref{c2}). The reference to Theorem \ref{Th5} gives the sufficiency of (\ref{d1}).

Let (\ref{c1}) hold. Then

$$\Bigl|\int_\Omega u\, \Bbb{V}_F(dx)\Bigr | ^2 \leq \Bbb{V}(F)\, \|\nabla u\|^2_{L_2(\Omega)}.$$

\noindent
Omitting a standard approximation argument, we put 
$$u(x) = \int_\Omega G(x,y)\, \Bbb{V}_F(dy)$$

\noindent
and the necessity of (\ref{d2}) results. $\qquad\square$

The next assertion follows directly from Theorem \ref{Th6}.

\begin{corollary}\label{Cor3}
The trace inequality (\ref{c5}) holds if and only if there exists a constant $C>0$ such that
\begin{equation}\label{k1}
\int_F\!\int_F G(x,y)\, \Bbb{V}_F(dx)\, \Bbb{V}_F(dy) \leq C\, \Bbb{V}(F)
\end{equation}
for all compact sets $F$ in $\Omega$.
\end{corollary}

\begin{remark}\label{3}
{\rm
 Note that Theorem \ref{Th6} and Corollary \ref{Cor3} include the case $n=2$ when the existence of the Green function is equivalent to ${\rm cap}(\Bbb{R}^2\backslash \Omega) >0$ (see \cite{[Lan]}).
}
\end{remark}

\begin{remark}\label{4}
{\rm
 Obviously, the pointwise estimate
\begin{equation}\label{d0}
\int_\Omega G(x,y)\, \Bbb{V}(dx) \leq 1/4
\end{equation}
implies (\ref{d1}) and hence it is sufficient for (\ref{c1}) to hold. $\qquad \square$
}
 \end{remark}

It is well-known that the operator $\tilde{S}$ obtained by the closure of the quadratic form $S[u,u]$  defined by (\ref{b0}) generates a contractive semigroup on $L_p(\Omega)$, $p\in (1,\infty)$ if and only if
\begin{equation}\label{d3b}
\int_\Omega |u|^{p-2}\, u\, \tilde{S}\, u\, dx \geq 0
\end{equation}
for all $u\in C^\infty_0(\Omega)$ (\cite{[LPh]}, \cite{[MS]}, \cite{[RS2]}, Th. X.48). The following analytic conditions related to (\ref{d3b}) can be deduced from Theorem \ref{Th4}.

\begin{corollary}\label{Cor7}
Let $p\in (1,\infty)$ and $p'= p/(p-1)$. The operator $\tilde{S}$ generates a contractive semigroup on $L_p(\Omega)$ if

$$\sup_{F\subset\Omega} \frac{\Bbb{V}(F)}{{\rm cap}_\Omega F} \leq \frac{1}{pp'}$$

\noindent
and only if
$$\sup_{F\subset\Omega} \frac{\Bbb{V}(F)}{{\rm cap}_\Omega F} \leq \frac{4}{pp'}.$$
\end{corollary}

\section{Trace inequality  for $\Omega = \Bbb{R}^n$}

Inequality
\begin{equation}\label{d3}
\int_{\Bbb{R}^n} |u|^2\, \Bbb{V}(dx) \leq C\, \int_{\Bbb{R}^n} |\nabla \, u|^2\, dx
\end{equation}
deserves to be discussed in more detail. First, 
 (\ref{d3}) for $n=2$ implies $\Bbb{V} =0$ . Let $n>2$. Needless to say, by Theorem \ref{Th4} the condition
\begin{equation}\label{d4}
\sup _F\frac{\Bbb{V}(F)}{{\rm cap}F} <\infty,
\end{equation}
where the supremum is taken over all compact sets $F$ in $\Bbb{R}^n$,
is necessary and sufficient  for (\ref{d3}). Restricting ourselves to arbitrary balls $B$ in $\Bbb{R}^n$, we have by (\ref{d4}) the obvious necessary condition
\begin{equation}\label{d5}
\sup _B\frac{\Bbb{V}(B)}{({\rm mes}_n B)^{1-2/n}} <\infty.
\end{equation}
On the other hand, using the isocapacitary inequality (\ref{a3}), we obtain the sufficient condition
\begin{equation}\label{d6}
\sup _F\frac{\Bbb{V}(F)}{({\rm mes}_n F)^{1-2/n}} <\infty,
\end{equation}
where the supremum is taken over all compact sets $F$ in $\Bbb{R}^n$. Moreover, the best value of $C$ in (\ref{d3}) satisfies
$$C\leq \frac{4v_n^{-2/n}}{n(n-2)} \sup _F\frac{\Bbb{V}(F)}{({\rm mes}_n F)^{1-2/n}}$$
and the constant factor in front of the supremum  is sharp \cite{[Maz3]}.

Although (\ref{d5}) and (\ref{d6}) look similar, they are not equivalent in general. In other words, one cannot replace arbitrary sets $F$ in (\ref{d4}) by balls. Paradoxically, the situation with the criterion (\ref{k1}) in the case $\Omega = \Bbb{R}^n$
is different. In fact,  Kerman and Sawyer \cite{[KeS]} showed that the trace inequality (\ref{d3}) holds if and only if for all balls $B$ in $\Bbb{R}^n$
\begin{equation}\label{k2}
\int_B\!\int_B\frac{\Bbb{V}(dx)\, \Bbb{V}(dy)}{|x-y|^{n-2}} \leq C\, \Bbb{V}(B).
\end{equation}

Maz'ya and Verbitsky  \cite{[MV1]} gave another necessary and sufficient condition for (\ref{d3}):
\begin{equation}\label{k3}
\sup_x\frac{I_1(I_1\Bbb{V})^2(x)}{I_1\, \Bbb{V}(x)} <\infty,
\end{equation}

\noindent
where $I_s$ is the Riesz potential of order $s$, i.e.
$$I_s\Bbb{V}(x): = \int_{\Bbb{R}^n} \frac{\Bbb{V}(dy)}{|x-y|^{n-s}}.$$

The following complete characterization of (\ref{d3}) is due to Verbitsky \cite{[Ver]}. For every dyadic cube $P_0$ in $\Bbb{R}^n$
\begin{equation}\label{k30}
\sup_{P_0}\frac{1}{V(P_0)}\sum \frac{\Bbb{V}(P)^2}{({\rm mes}_n P)^{1-2/n}} <\infty,
\end{equation}
where the sum is taken over all dyadic cubes $P$ contained in $P_0$ and $C$ does not depend on $P_0$.

Let us assume that $\Bbb{V}$ is absolutely continuous with respect to ${\rm mes}_n$, i.e.
$$\Bbb{V}(F) = \int_F V(x)\, dx.$$
 Fefferman and Phong \cite{[F]} proved that  (\ref{d3}) follows from
\begin{equation}\label{e0}
\sup_B\frac{\displaystyle{\int_B V^t\, dx}}{({\rm mes}_n B)^{1-2t/n}} < \infty \qquad {\rm with}\,\, {\rm some}\,\,\,
 t\in (1, n/2).
\end{equation}
This can be readily deduced from Sawyer's inequality  \cite{[Saw]}
\begin{equation}\label{e1}
\int_{\Bbb{R}^n} ({\cal M}_t f)^2\, \nu(dx) \leq C\sup_B\frac{\nu(B)}{({\rm mes}_n B)^{1-2t/n}}\int_{\Bbb{R}^n} f^2\, dx,
\end{equation}
where $\nu$ is a measure,
$f\geq 0$ and ${\cal M}_t $ is the fractional maximal operator defined by
$${\cal M}_t f(x) = \sup_{B} \frac{\displaystyle{\int_B f(y)\, dy}}{({\rm mes}_n B)^{1-t/n}}.$$
 In fact,  for any $\delta>0$
$$I_1f(x) = (n-1)\int_0^\infty r^{-n}\int_{B_r(x)} f(y)\, dy\,
dr\leq (n-1)\bigl(\delta\, {\cal M} f(x) + (t-1)^{-1}\delta^{1-t}
{\cal M}_t f(x)\bigr),$$ where ${\cal M}$ is the Hardy-Littlewood
maximal operator. Minimizing the right-hand side in $\delta$, we
arrive at the  inequality \begin{equation}\label{e2}
I_1 f(x) \leq \frac{(n-1)t}{t-1}\, ({\cal M}_t f(x))^{1/t} ({\cal
M}\, f(x))^{1-1/t}
\end{equation}
(\cite{[Hed]}, see also \cite{[AH]}, Sect.\,3.1).  Now, (\ref{e2}) implies
$$\|I_1\, f\|_{L_2(Vdx)} \leq C\, \|{\cal M}_t f\|^{1/t}_{L_2(V^tdx)} \|{\cal M}\, f\|^{1-1/t}_{L_2(\Bbb{R}^n)}.$$
Therefore, by (\ref{e1}) and by the boundedness of ${\cal M}$ in $L_2(\Bbb{R}^n)$ we arrive at the  inequality
$$\|I_1\, f\|_{L_2(Vdx)} \leq C\, \sup_B\Biggl(\frac{\displaystyle{\int_B V^t\, dx}}{({\rm mes}_n B)^{1-2t/n}}\Biggr)^{1/2t} \, \|f\|_{L_2(\Bbb{R}^n)}$$
which is equivalent to the Fefferman-Phong result. Their result was improved by Chang, Wilson and Wolff  \cite{[ChWW]} who showed that if $\varphi$ is an increasing function $[0,\infty) \to [1,\infty)$ subject to
\begin{equation}\label{e3}
\int_1^\infty (\tau\, \varphi(\tau))^{-1}d\tau <\infty,
\end{equation}
then the condition
\begin{equation}\label{e4}
\sup_B\frac{\displaystyle{\int_B V(x)\, \varphi\bigl(V(x)({\rm diam}\, B)^2\bigr)\, dx}}{({\rm mes}_n\, B)^{1-2/n}} <\infty
\end{equation}
is sufficient for (\ref{d3}) to hold. The assumption (\ref{e3}) is sharp.

By Remark \ref{4}, the condition
\begin{equation}\label{e00}
\sup_x\int_{\Bbb{R}^n}\frac{\Bbb{V}(dy)}{|x-y|^{n-2}} <\infty,
\end{equation}
is sufficient for (\ref{d3}). It is related to Kato's condition frequently used in the spectral theory of the Schr\"odinger operator
(see Sect.\,14  below for more information about Kato's condition).  Note that (\ref{e00}) is rather far from being necessary since it excludes (\ref{d3}) with  $\Bbb{V}(dx) = |x|^{-2} dx$, i.e. the classical Hardy inequality, unlike other  conditions just mentioned.

We conclude this section with the observation that the multiplicative inequality
\begin{equation}\label{x1}
\int_{\Bbb{R}^n} u^2\, \Bbb{V}(dx) \leq C\, \Bigl(\int_{\Bbb{R}^n}|\nabla \, u|^2 dx\Bigr)^\tau\Bigl(\int_{\Bbb{R}^n} u^2\, dx\Bigr)^{1-\tau}, \qquad 0\leq \tau <1,
\end{equation}
 is equivalent to
\begin{equation}\label{x2}
\sup_B\frac{\Bbb{V}B}{({\rm mes}_n B)^{1-2\tau/n}} <\infty
\end{equation}
(\cite{[Maz7]}, Theorem 1.4.7, see Sect. 15 below for
further development).

\section{Positive solutions of $(\Delta + \Bbb{V})w =0$}

Inequality (\ref{d3}) is intimately related to the problem of existence of positive solutions for the Schr\"odinger equation
\begin{equation}\label{x2a}
-\Delta\, w = \Bbb{V}\, w \qquad {\rm in}\,\,\, \Bbb{R}^n.
\end{equation}
 S.\,Agmon showed that the existence of $w\geq 0$ is equivalent to the positivity of the operator $-\Delta -\Bbb{V}$ for relatively nice $\Bbb{V}$ \cite{[Agm]}. The following result due to Hansson, Maz'ya and Verbitsky \cite{[HMV]} shows that the condition
\begin{equation}\label{x2b}
I_1(I_1\Bbb{V})^2 \leq C\, I_1\, \Bbb{V}
\end{equation}
is equivalent (up to values of $C$)  to the existence of $w>0$.

\begin{theorem}\label{Th6}
{\rm (i)} If $-\Delta\, w = \Bbb{V}\, w$ has a nonnegative (weak) solution $w$, then $I_1\Bbb{V} <\infty$ almost everywhere and there exists a constant $C_1 = C(n)$ such that
\begin{equation}\label{e6}
I_1\bigl( (I_1\Bbb{V})^2\bigr) (x) \leq C_1\, I_1\, \Bbb{V}(x) \quad {\rm a.e.}
\end{equation}

{\rm (ii)} Conversely, there exists a constant $C_2 = C_2(n)$ such that if {\rm(\ref{e6})} holds with $C_2$ in place of $C_1$, then there exists a positive solution $w$ to {\rm (\ref{x2a})} which satisfies
$$|\nabla\, \log\, w(x)| \leq C\, I_1\, w(x)$$
and
$$w(x) \geq C\, (|x| +1)^{-c}$$
for some positive constants $C$ and $c$.
\end{theorem}

In addition, if $I_2\Bbb{V}(x) <\infty$ a.e., then there is a solution $w$ such that
$$e^{I_2\, \Bbb{V}(x)} \leq w(x) \leq e^{C\, I_2\, \Bbb{V}(x)}.$$
All the constants depend only on $n$.

\section{Semiboundedness of the operator $-\Delta -\Bbb{V}$}

We formulate some consequences of Theorem \ref{Th4} concerning the
topic in the title, i.e. the inequality

$$\int_\Omega |\nabla \, u|^2\, dx - \int_\Omega |u|^2\, \Bbb{V}(dx) \geq -C\, \int_\Omega |u|^2\, dx, \quad u\in C^\infty_0(\Omega).$$

\begin{theorem}\label{Th7}
{\rm (\cite{[Maz3]}, see also \cite{[Maz7]},
Sect. 2.2)}

{\rm (i)} If
\begin{equation}\label{k4}
\mathop{\hbox {{\rm lim sup}}}_{\delta \to 0}\, \Bigl\{\frac{\Bbb{V}(F)}{{\rm cap}_\Omega F}: F\subset\Omega, \,\, {\rm diam}\, F \leq \delta\Bigr\} <\frac{1}{4},
\end{equation}
then the quadratic form $S[u,u]$ defined by {\rm (\ref{b0})} is semi-bounded from below and closable in $L_2(\Omega)$.

{\rm (ii)} If the form $S[u,u]$ is semi-bounded from below in $L_2(\Omega)$, then

\begin{equation}\label{k5}
\mathop{\hbox {{\rm lim sup}}}_{\delta \to 0}\Bigl\{\frac{\Bbb{V}(F)}{{\rm cap}_\Omega F}: F\subset\Omega, \,\, {\rm diam}\, F \leq \delta\Bigr\} \leq 1.
\end{equation}
\end{theorem}

\begin{corollary}\label{cor1}
The condition
\begin{equation}\label{k6}
\mathop{\hbox {{\rm lim sup}}}_{\delta \to 0}\, \Bigl\{\frac{\Bbb{V}(F)}{{\rm cap}_\Omega F}: F\subset\Omega, \,\, {\rm diam}\, F \leq \delta\Bigr\}  =0
\end{equation}
is necessary and sufficient for the semiboundedness of the form
\begin{equation}\label{f0}
S_h[u,u]: = h\int_\Omega |\nabla \, u|^2 dx - \int_\Omega |u|^2\, \Bbb{V}(dx)
\end{equation}
in $L_2(\Omega)$ for all $h>0$.
\end{corollary}

\begin{corollary}\label{cor2}
The inequality
\begin{equation}\label{k7}
\int_\Omega |u|^2\, \Bbb{V}(dx) \leq C\, \int_\Omega(|\nabla \, u|^2 + |u|^2)\, dx,
\end{equation}
where $u$ is an arbitrary function in $C^\infty_0(\Omega)$ and $C$ is a constant independent of $u$,  holds if and only if there exists $\delta>0$ such that
$$\sup\Bigl \{\frac{\Bbb{V}(F)}{{\rm cap}_\Omega F}: \, F\subset\Omega, \,\, {\rm diam}\, F\leq \delta\Bigr\} <\infty.$$
\end{corollary}

\section{Negative spectrum of $-h\Delta - \Bbb{V}$}

Investigation of the negative spectrum of the Schr\"odinger operator with negative potential can be based upon the following two classical general results.

\begin{lemma}\label{lem1}
{\rm  \cite{[Fr]}} Let $A[u,u]$ be a closed symmetric quadratic
form in a Hilbert space $H$ with domain $D[A]$ and let $\gamma(A)$
be its positive greatest lower bound. Further, let $B[u,u]$ be
a real  valued  quadratic form, compact in $D[A]$. Then the form
$A-B$ is semi-bounded below in $H$, closed in $D[A]$, and its
spectrum is discrete to the left of $\gamma(A)$.
\end{lemma}

\begin{lemma}\label{lem2}
{\rm  \cite{[Gl]}} For the negative spectrum of a selfadjoint
operator $A$ to be infinite it is necessary and sufficient that
there exists a linear manifold of infinite dimension on which
$(A\, u, u)<0$.
\end{lemma}

Note also that by Birman's general theorem  \cite{[Bir2]}, the
discreteness  of the negative spectrum of the operator
$$\tilde{S}_h = -h\Delta - \Bbb{V} \qquad {\rm for}\,\,{\rm all}\,\, h>0$$
generated by the closure of the quadratic form $S_h[u,u] $ defined by (\ref{f0})
is equivalent to the compactness of the quadratic form
\begin{equation}\label{k9a}
\int_\Omega |u|^2\, \Bbb{V}(dx)
\end{equation}
with respect to the norm
 $$\Bigl(\int_\Omega(|\nabla \, u|^2 + |u|^2)\, dx\Bigr)^{1/2}.$$
Analogously, the finiteness of the negative spectrum of the operator $\tilde{S}_h$ is equivalent to the compactness of the quadratic form (\ref{k9a}) with respect to the norm $\|\nabla \, u\|_{L_2(\Omega)}$.

Moreover, by \cite{[Bir1]} (see also \cite{[Bir2]} and
\cite{[Schw]}) the number of negative eigenvalues of $\tilde{S}_h$
coincides with the number of eigenvalues $\lambda_k$, $\lambda_k <
h^{-1}$, of the Dirichlet problem
$$ -\Delta u -
\lambda \Bbb{V} \, u = 0, \quad u\in \ring L_2^1(\Omega).$$

The next theorem contains analytic conditions both necessary and sufficient for the Schr\"odinger operators $\tilde{S}_h$ to have discrete, infinite or finite negative spectra for all $h>0$.
These characterizations were obtained by the author   from either sufficient or necessary conditions, analogous to (\ref{k4}), (\ref{k5}), for the operator $\tilde{S}$ independent of the parameter $h$ to have a negative spectrum with the properties just listed (\cite{[Maz1]}, \cite{[Maz3]}, see also  \cite{[Maz7]}, Sect. 2.2.)

\begin{theorem}\label{Th8}
 {\rm (i)} A necessary and sufficient condition for the discreteness of the negative spectrum of $\tilde{S}_h$  for all $h>0$ is
\begin{equation}\label{k9}
\sup_{F\subset\Omega\backslash B_\rho, \, {\rm diam} F\leq 1} \frac{\Bbb{V}(F)}{{\rm cap}_\Omega F} \to 0 \quad {\rm as}\,\, \rho\to \infty.
\end{equation}
Here and elsewhere $B_\rho = \{x\in \Bbb{R}^n: |x| <\rho\}$.

{\rm (ii)} A necessary and sufficient condition for the discreteness of the negative spectrum of $\tilde{S}_h$  for all $h>0$ is

\begin{equation}\label{k10}
\sup_{F\subset\Omega\backslash B_\rho} \frac{\Bbb{V}(F)}{{\rm cap}_\Omega F} \to 0 \quad {\rm as}\,\, \rho\to \infty.
\end{equation}

{\rm (ii)} A necessary and sufficient condition for the negative spectrum of $\tilde{S}_h$ to be infinite for all $h>0$ is

\begin{equation}\label{k11}
\sup_{F\subset\Omega} \frac{\Bbb{V}(F)}{{\rm cap}_\Omega F} = \infty.
\end{equation}
\end{theorem}

Needless to say, simpler sufficient conditions with ${\rm cap}_\Omega F$ replaced by a function of ${\rm mes}_n F$ follow directly from (\ref{k9}) - (\ref{k11}) and the isocapacitary inequalities (\ref{a1}) - (\ref{a3}).

\section{Rellich-Kato  theorem}

By the basic Rellich-Kato  theorem \cite{[Re1]} - \cite{[Ka1]}, the self-adjointness of $-\Delta + V$ in $L_2(\Bbb{R}^n)$ is guaranteed by the inequality
\begin{equation}\label{g1}
\|V\, u\|_{L_2(\Bbb{R}^n)} \leq a\, \|\Delta\, u\|_{L_2(\Bbb{R}^n)} + b\, \|u\|_{L_2(\Bbb{R}^n)},
\end{equation}
where $a<1$ and $u$ is an arbitrary function in $C^\infty_0 (\Bbb{R}^n)$.

Let $n\leq 3$. The Sobolev embedding $W^{2}_2(\Bbb{R}^n) \subset L_\infty(\Bbb{R}^n)$ and an appropriate choice of the test function in (\ref{g1}) show that (\ref{g1}) holds if and only if there is a sufficiently small constant $c(n)$ such that
\begin{equation}\label{k31}
\sup_{x\in \Bbb{R}^n} \int_{B_1(x)} |V(y)|^2\, dy \leq c(n).
\end{equation}
Here and elsewhere $B_r(x) = \{y\in \Bbb{R}^n:|y-x|<r\}$.  This
class of potentials is known as Stummel's class $S_n$ which is
defined also for higher dimensions by
$$\lim_{r\downarrow 0} \Bigl(\sup_x\int_{B_r(x)} |x-y|^{4-n}|V(y)|^2\, dy\Bigr) =0\qquad\, \,  {\rm for}\,\, n\geq 5,$$
$$\lim_{r\downarrow 0} \Bigl(\sup_x\int_{B_r(x)} \log(|x-y|^{-1})\, |V(y)|^2\, dy\Bigr) =0\qquad {\rm for}\,\, n=4,$$
\cite{[Stum]}. Although the condition $V\in S_n$ is sufficient for
(\ref{g1}) for every $n$,  it does not seem quite natural for
$n\geq 4$. As a matter of fact, it excludes the simple potential
$V(x) = c|x|^{-2}$ 
obviously satisfying (\ref{g1}), if the factor  $c$ is small enough.

If $n\geq 5$, a characterization of (\ref{g1}) (modulo best constants) results directly from a necessary and sufficient condition for the inequality
$$\|V\, u\|_{L_2(\Bbb{R}^n)} \leq C\, \|\Delta \, u\|_{L_2(\Bbb{R}^n)}, \qquad u\in C^\infty_0(\Bbb{R}^n),$$
found in \cite{[Maz4]}. We claim that the Rellich-Kato condition (\ref{g1}) holds in the case $n\geq 4$ if and only if there is a sufficiently small $c(n)$ subject to
\begin{equation}\label{g2}
\sup_{{\rm diam} F\leq 1} \frac{\displaystyle{\int_F|V(y)|^2\, dy}}{{\rm cap}_2 F} \leq c(n)
\end{equation}
(the values of $c(n)$ in the sufficiency and necessity parts are different similarly to the criteria in Sect. 4).
 Here  and elsewhere ${\rm cap}_m$ is the polyharmonic capacity of a compact set $F$ defined in the case $2m<n$ by
  $${\rm cap}_m F = \inf \Bigl\{\int_{\Bbb{R}^n} |\nabla_m\, u|^2\, dx:\,\, u\in C^\infty_0(\Bbb{R}^n), \, u\geq 1 \, {\rm on}\,  F\Bigr\},$$
   where $\nabla _m = \{\partial^m/\partial x_1^{\alpha_1}\ldots \partial x_n^{\alpha_n}\}$. For $2m =n$  one should add $\|u\|_{L_2(\Bbb{R}^n)}^2$ to the last integral.

 The condition
\begin{equation}\label{g2a}
\sup_{{\rm diam} F\leq \delta} \frac{\displaystyle{\int_F|V(y)|^2\, dy}}{{\rm cap}_2 F} \to 0 \quad {\rm as}\,\, \delta\to 0
\end{equation}
is necessary and sufficient for (\ref{g1})  to hold with an arbitrary $a$ and $b = b(a)$.

 An obvious necessary condition for (\ref{g1}) is
\begin{equation}\label{g3}
\cases{\displaystyle{ \sup_{r\leq 1, x\in \Bbb{R}^n}
r^{4-n}\int_{B_r(x)} |V(y)|^2\, dy \leq c(n) \qquad\qquad {\rm
for}}\,\, n\geq 5\nonumber\cr \nonumber\cr
\displaystyle{\sup_{r\leq 1, x\in \Bbb{R}^n}
\bigl(\log\frac{2}{r}\bigr)^{-1}\int_{B_r(x)} |V(y)|^2\, dy \leq
c(n) \qquad {\rm for}\,\, n=4}}
\end{equation}
where $c(n)$ is sufficiently small.
Standard lower estimates of ${\rm cap}_2$ by ${\rm mes}_n$ combined with the criterion (\ref{g2}) give the sufficient condition
\begin{equation}\label{g4}
\cases{\displaystyle{
\sup_{ {\rm diam} F \leq 1} ({\rm mes}_n F)^{(4-n)/n}\int_{F} |V(y)|^2\, dy \leq c(n) \qquad {\rm for}}\,\, n\geq 5\nonumber\cr
\nonumber\cr
\displaystyle{\sup_{ {\rm diam} F \leq 1} \bigl(\log\frac{v_n}{{\rm mes}_n F}\bigr)^{-1}\int_{F} |V(y)|^2\, dy \leq c(n) \qquad\,  {\rm for}\,\, n=4.}}
\end{equation}
Though sharp and looking  similar, (\ref{g3}) and (\ref{g4}) are not equivalent. We omit a discussion of other sufficient and both necessary and sufficient conditions for  (\ref{g2}) to hold parallel to that in Sect.\,5 (analogs of the Fefferman $\&$ Phong, Kerman $\&$ Sawyer,  Maz'ya $\&$ Verbitsky, and Verbitsky conditions mentioned in Sect. 5).

Finally, we observe that the inequality
$$\|V\, u\|_{L_2(\Bbb{R}^n)} \leq C\, \|\Delta\, u\|^\tau_{L_2(\Bbb{R}^n)}\|u\|^{1-\tau}_{L_2(\Bbb{R}^n)}$$
holds for a certain $\tau \in (0,1)$ and every $u\in C^\infty_0(\Bbb{R})$ if and only if for all $r\in (0,1)$
$$\sup_x\int_{B_r(x)} |V(y)|^2\, dy \leq C\, r^{n-4\tau}$$
(see Theorem 1.4.7 in \cite{[Maz7]}).

\section{Sobolev regularity for solutions of $(-\Delta + V)u = f$}

Let $m$ be integer $\geq 2$ and let $W_2^m(\Bbb{R}^n)$ denote the Sobolev space of functions $u\in L_2(\Bbb{R}^n)$ such that $\nabla _m u\in L_2(\Bbb{R}^n)$.
Obviously, the operator
$$-\Delta +V : W^{m}_2(\Bbb{R}^n) \to W_2^{m-2}(\Bbb{R}^n)$$
is bounded if and only if $V$ is a pointwise multiplier acting
from $W^{m}_2(\Bbb{R}^n)$ into $W_2^{m-2}(\Bbb{R}^n)$. (We use the
notation $V\in M(W^{m}_2(\Bbb{R}^n) \to W^{m-2}_2(\Bbb{R}^n)$).
According to  \cite{[MSha1]}, Ch.\,1,   necessary and  sufficient
conditions  for $V\in M(W^{m}_2(\Bbb{R}^n) \to
W^{m-2}_2(\Bbb{R}^n)$ have the form
\begin{equation}\label{m1}
\sup_{x\in \Bbb{R}^n} \int_{B_1(x)}(|\nabla_{m-2} V(y)|^2+
|V(y)|^2)\, dy <\infty
\end{equation}
for $n<2m$ and
\begin{equation}\label{m2}
\sup_{{\rm diam} F\leq 1}
\frac{\displaystyle{\int_F|\nabla_{m-2}V(y)|^2\, dy}}{{\rm cap}_m
F} +\sup_{x\in \Bbb{R}^n} \int_{B_1(x)} | V(y)|^2\, dy <\infty
\end{equation}
for $n\geq 2m$.

Two-sided estimates for the essential norm of a multiplier:
$W^{m}_2(\Bbb{R}^n) \to W^{m-2}_2(\Bbb{R}^n)$ found in
\cite{[MSha1]}, Ch.\,4,  lead to the following regularity result.

\begin{theorem}\label{Thm1}
Let $\varphi\in W^{1}_2(\Bbb{R}^n, loc)$ be a solution of
$$(-\Delta + V)\varphi = f,$$
where $f\in W^{m-2}_2(\Bbb{R}^n, loc)$.

If $n<2m$ and $V\in W^{m-2}_2(\Bbb{R}^n,loc)$ then $\varphi \in W^{2m}_2(\Bbb{R}^n, loc)$.

Let $n\geq 2m$. Suppose   there exists a sufficiently small constant $c(n)$ such that for all sufficiently small $\delta>0$
\begin{equation}\label{m3}
 \sup_{{\rm diam} F \leq \delta}\frac{\displaystyle{\int_F|\nabla_{m-2}V(y)|^2\, dy}}{{\rm cap}_m F} +
\sup_x\frac{\displaystyle{\int_{B_\delta(x)}|V(y)|^2\, dy}}{\delta^{n-2m}}
\leq c(n). 
\end{equation}
Then $\varphi \in W^{2m}_2(\Bbb{R}^n, loc)$.
\end{theorem}

Equivalent forms of the conditions (\ref{m2}) and (\ref{m3})  of Kerman $\&$ Sawyer, Maz'ya $\&$ Verbitsky, and Verbitsky type,  as well as various separately necessary or sufficient  conditions involving no capacities,  are available (compare with  Sect. 5).

 Similar higher regularity $L_p$-results can be obtained  from estimates of the essential norms of multipliers:
 $W^{m}_p(\Bbb{R}^n) \to W^{m-2}_p(\Bbb{R}^n)$ obtained in \cite{[MSha1]}, Ch. 4.

\section{Relative form boundedness and form compactness}

 Maz'ya and Verbitsky \cite{[MV2]}  gave
necessary and
 sufficient  conditions for the  inequality
\begin{equation}\label{12.1}
\left \vert \int_{\Bbb{R}^n} |u(x)|^2 \, V(x) \, dx\right
\vert \leq C
\int_{\Bbb{R}^n} |\nabla u(x)|^2  \, dx, \quad u
\in C_0^\infty (\Bbb{R}^n)
\end{equation}
to hold. Here
 the ``indefinite weight'' $V$ may change sign,
or even be a complex-valued
 distribution on $\Bbb{R}^n$, $n \ge 3$. (In the  latter case, the
left-hand
side of (\ref{12.1})
is understood as $|<V u, u>|$, where $<V \, \cdot, \, \cdot>$ is the
quadratic
form associated with the corresponding multiplication operator $V$.)  An analogous inequality   for the
Sobolev space $W_2^1 (\Bbb{R}^n)$,  $n \ge 1$ was also characterized in \cite{[MV2]}:
\begin{equation}\label{12.2}
\left \vert \int_{\Bbb{R}^n} |u(x)|^2 \, V(x) \, dx\right
\vert
 \leq C
\int_{\Bbb{R}^n} \bigl( |\nabla u(x)|^2 + |u(x)|^2 \bigr)  \, dx, \quad u
\in C_0^\infty (\Bbb{R}^n).
\end{equation}

Such inequalities  are
used extensively in  spectral and scattering
theory
of the Schr\"odinger operator $H_V = -\Delta + V$ 
and its higher-order analogs, especially in questions of
self-adjointness, resolvent convergence, estimates for the number
of bound states,  Schr\"odinger semigroups,
 etc. (See \cite{[Bir2]},   \cite{[BS1]},  \cite{[BS2]}, \cite{[CZh]},
\cite{[Dav1]}, \cite{[Far]},  \cite{[F]}, \cite{[RS2]},
\cite{[Sch1]}, \cite{[Sim3]}, and the literature cited there.) In
particular,  (\ref{12.2})
 is equivalent
 to the fundamental  concept
of the relative boundedness of the
 potential energy operator $V$ with respect to $ - \Delta$
in the sense of quadratic forms. Its abstract version
appears in the so-called KLMN Theorem,  which is
discussed in detail together with  applications to quantum-mechanical
Hamiltonian operators, in \cite{[RS2]}, Sec. X.2.

It follows from the polarization identity that  (\ref{12.1})
can be restated equivalently in terms of the  corresponding
sesquilinear form:
$$|<V u, v>| \leq C  \,
||\nabla u||_{L_2} \,  \,
||\nabla v||_{L_2},
$$
for all $u, v \in C_0^\infty (\Bbb{R}^n)$.  In other words, it
is equivalent to the
boundedness
of the operator $H_V$,
\begin{equation}\label{12.3}
H_V : \, \ring L^1_2 (\Bbb{R}^n) \to L_2^{-1}(\Bbb{R}^n), \qquad n \ge 3.
\end{equation}
Here  the energy space
$\ring L^1_2 (\Bbb{R}^n)$ is defined as the completion of $C_0^\infty (\Bbb{R}^n)$ with
respect
to the Dirichlet norm $||\nabla u||_{L_2 }$, and $L_2^{-1}(\Bbb{R}^n)$
is the dual of $\ring L^1_2(\Bbb{R}^n)$.
Similarly, (\ref{12.2}) means that
$H$ is a bounded operator which maps
$W_2^1 (\Bbb{R}^n)$ to $W_2^{-1}(\Bbb{R}^n), \, n \ge 1$.

Previously, the case
of {\it nonnegative\/}  $V$ in (\ref{12.1}) and (\ref{12.2})
has been studied in a comprehensive way (see Sect.\,5). For general
$V$, only sufficient conditions have been known.

It is worthwhile to
 observe that the usual ``na\"\i ve'' approach
 is to decompose
 $V$
into its positive and negative parts: $V = V_+ - V_-$, and to apply the
just
mentioned results to both $V_+ $ and $V_- $. However, this procedure
drastically
 diminishes the class of admissible weights  $V$ by ignoring
a possible cancellation between $V_+ $ and $V_- $.
This cancellation phenomenon is evident for
 strongly oscillating weights considered below. Examples of this
type are
known mostly in relation
to quantum mechanics problems  \cite{[ASi]}, \cite{[CoG]},
\cite{[NS]},  \cite{[Stu]}.

 The following result obtained in \cite{[MV2]} reflects
 a general principle
 which has  much wider range of applications.

 \begin{theorem}\label{Theorem 12.1}
 Let $V$ be a complex-valued
 distribution on $\Bbb{R}^n$, $n \ge 3$. Then
  {\rm (\ref{12.1})} holds  if and only if $V$
is the divergence
of a vector-field $\vec \Gamma : \, \Bbb{R}^n \to \Bbb{C}^n$  such that
\begin{equation}\label{12.4}
\int_{\Bbb{R}^n} |u(x)|^2 \, |\vec \Gamma (x)  |^2  \, dx
\leq {\rm  {const}}
\int_{\Bbb{R}^n} |\nabla u(x)|^2  \, dx,
\end{equation}
where the constant is independent of
 $u \in C_0^\infty (\Bbb{R}^n)$. The vector-field
$\vec \Gamma  \in \bold L_{2}(\Bbb{R}^n,loc)$ can be
chosen as  $\vec \Gamma  = \nabla \Delta^{-1} V$.

Equivalently,  the
Schr\"odinger operator $H_V$ acting from
$\ring L^1_2 (\Bbb{R}^n)$ to $L_2^{-1}(\Bbb{R}^n)$ is bounded if and only if
{\rm (\ref{12.4})} holds.
Furthermore, the corresponding multiplication operator $V: \,
\ring L^1_2 (\Bbb{R}^n) \to L_2^{-1}(\Bbb{R}^n)$ is  compact if and only if the
embedding
$$\ring L^1_2 (\Bbb{R}^n)\,  \subset \, L_2( \Bbb{R}^n, \, |\vec \Gamma|^2 \, dx)
$$
is compact.
\end{theorem}

Once $V$ is written
as  $V = {\rm  div} \, \vec \Gamma$,
the implication (\ref{12.4})$\rightarrow $(\ref{12.1}) becomes trivial. In fact,
it follows
 using integration by parts and the Schwarz inequality.
This idea has been known for a long time in mathematical physics
(see, e.g.,   \cite{[CoG]}) and
  theory of multipliers in Sobolev spaces
\cite{[MSha1]}.  On the other
hand, the proof of the converse statement (\ref{12.1})$\rightarrow$(\ref{12.4})
where $\vec \Gamma  = \nabla \Delta^{-1} V$
 is  rather delicate.

Theorem \ref{Theorem 12.1} makes it possible to reduce
the problems of boundedness and compactness
for general ``indefinite'' $V$
 to the case of  nonnegative  weights
$|\vec \Gamma|^2$,
which is by now
well understood.  In particular,
 combining  Theorem \ref{Theorem 12.1} and the  criteria in Sect.\,5 for $|\vec \Gamma|^2$ one arrives at  analytic necessary and sufficient conditions for (\ref{12.1}) to hold.

As a corollary, one obtains a  necessary condition
for (\ref{12.1})  in terms of Morrey spaces of negative order.

\begin{corollary}\label{Corollary 12.1} If {\rm  (\ref{12.1})} holds,  then,
 for every  ball $B_r(x_0)$ of radius $r$,
$$
\int_{B_r(x_0)} |\nabla \Delta^{-1} V (x)|^{2} \, dx
\leq C \, r^{n - 2 },
$$
where the constant
does not depend on $x_0 \in \Bbb{R}^n$ and $r>0$.
\end{corollary}

\begin{corollary}\label{Corollary 12.2}
 In the statements of Theorem \ref{Theorem 12.1}  and Corollary \ref{Corollary 12.1},
 one can put
the scalar function $(- \Delta)^{-\frac 1 2} V$ in place of
 $\vec \Gamma
 = \nabla \Delta^{-1} V$.
In particular,
 {\rm  (\ref{12.4})} is equivalent to the inequality:
\begin{equation}\label{12.5}
\int_{\Bbb{R}^n} |u(x)|^2 \, |(- \Delta)^{-\frac 1 2} V (x)  |^2  \, dx
\leq C
\int_{\Bbb{R}^n} |\nabla u(x)|^2  \, dx,
\end{equation}
for all $u \in C_0^\infty (\Bbb{R}^n)$.
\end{corollary}
The proof of Corollary \ref{Corollary 12.2} uses the boundedness
of  standard singular integral operators in the space of
functions $f \in L_{2,{{\rm loc}}}(\Bbb{R}^n)$ such that
$$
\int_{\Bbb{R}^n} |u(x)|^2 \, |f (x)|^2  \, dx
\leq C
\int_{\Bbb{R}^n} |\nabla u(x)|^2  \, dx,
$$
for all $u \in C_0^\infty (\Bbb{R}^n)$; this fact
 was established earlier in \cite{[MV1]}.

Corollary \ref{Corollary 12.2} indicates that
 a  decomposition into a positive and negative part, appropriate 
for (\ref{12.1}), should involve expressions like
$(-\Delta)^{-\frac 1 2} V$ rather than $V$
itself.
Another important consequence   is
that the class of weights $V$ satisfying (\ref{12.1}) is invariant under
 standard singular integral
and maximal operators.
\begin{remark}\label{Remark 1}
\rm{
Similar results are valid for inequality
(\ref{12.2}); one only has to replace
the operator
 $(- \Delta)^{-1/2}$ with   $ (1 - \Delta)^{-1/2}$,
and the Wiener capacity
${\rm {cap}} \, F$ with the corresponding Bessel capacity (see \cite{[AH]}).
 It suffices to restrict oneself, in the corresponding
statements,
 to cubes  or balls whose volumes are
less than $1$.
}
\end{remark}

We now discuss some related results  in terms of
more conventional
classes of admissible weights $V$.
The following corollary, which is
an immediate consequence of Theorem \ref{Theorem 12.1} and Corollary \ref{Corollary 12.2},
  gives a simpler sufficient
condition for
(\ref{12.1})  with  Lorentz-Sobolev spaces of negative order in its formulation.

\begin{corollary}\label{Corollary 12.3}
Let $n \ge 3$, and let $V$ be a
distribution
 on $\Bbb{R}^n$ such that
$(- \Delta)^{-\frac 1 2} V  \in L_{n, \, \infty} (\Bbb{R}^n)$,
where $L_{p, \, \infty}$ denotes the usual Lorentz
(weak $L_p$) space.
 Then {\rm  (\ref{12.1})}
holds.
\end{corollary}

For the definition and basic properties of Lorentz spaces
$L_{p,q} (\Bbb{R}^n)$ we
refer  to \cite{[SW]}. Note that, in particular,
$(- \Delta)^{-\frac 1 2} V
 \in L_{n, \, \infty}$ is equivalent to the
estimate
\begin{equation}\label{12.9}
\int_F |(- \Delta)^{-\frac 1 2} V (x)|^2 \, dx \le
C \,
 ({\rm mes}_n F)^{1-  2/n}.
 \end{equation}

The following corollary of Theorem \ref{Theorem 12.1} is applicable to
distributions $V$,
and encompasses a   class
of weights  which is broader than the Fefferman$-$Phong class (\ref{e0}) even
in the case where $V$ is a nonnegative measurable function.

\begin{corollary}\label{Corollary 12.4}
Let $V$ be a distribution on $\Bbb{R}^n$
which satisfies,
 for some $t > 1$, the inequality
\begin{equation}\label{12.11}
\int_{B_r(x_0)} |(- \Delta)^{-\frac 1 2} V (x)  |^{2t} \, dx
\leq C \, r^{n - 2 t},
\end{equation}
for every ball $B_r(x_0)$ in $\Bbb{R}^n$.
Then {\rm  (\ref{12.1})} holds.
\end{corollary}

 Note that by Corollary \ref{Corollary 12.1} the preceding
inequality with $t=1$
is necessary for  (\ref{12.1}) to hold.

\begin{remark}\label{Remark 4}
\rm{
 A refinement
of (\ref{12.11}) in terms of the  condition (\ref{e4})
established by Chang, Wilson, and Wolff  \cite{[ChWW]}
 is readily
 available by combining (\ref{e4}) with  Theorem \ref{Theorem 12.1}.
}
\end{remark}

To clarify the  multi-dimensional characterizations
for ``indefinite weights'' $V$ presented above,
 we state an elementary  analog of
Theorem \ref{Theorem 12.1}  for the
Sturm-Liouville operator  $H_V = - \frac{d^2} {d x^2} + V$
on the half-line.

\begin{theorem}\label{Theorem 12.3}
The inequality
\begin{equation}\label{12.12}
\left \vert \int_{\Bbb{R}_+} |u(x)|^2 \, V(x) \, dx\right
\vert \leq C
\int_{\Bbb{R}_+} |u' (x)|^2  \, dx,
\end{equation}
holds for all $u \in C^\infty_0 ( \Bbb{R}_+ ) $ if and only if
\begin{equation}\label{12.13}
\sup_{a>0} \, \, a \, \int_a^\infty
\left |\int_x^\infty V(t) \, dt\right |^2 \, dx  < \infty,
\end{equation}
where $\Gamma (x)  = \displaystyle{\int_x^\infty} V(t) \, dt$
is understood in terms of distributions.

Equivalently,  $H_V: \,
\ring L^1_2 (\Bbb{R}_+) \to L_2^{-1}(\Bbb{R}_+)$ is bounded if and only if
{\rm  (\ref{12.13})} holds.
Moreover, the corresponding multiplication operator
$V$
is  compact
 if and only if
\begin{equation}\label{12.14}
 a \, \int_a^\infty |\Gamma (x)|^2 \, dx   = o \, (1), \quad
 \,\,  {\rm  {where}} \quad a \to 0^+ \quad  {\rm  {and}}
\quad a \to +
 \infty.
\end{equation}
\end{theorem}
For nonnegative $V$,  condition (\ref{12.13}) is easily seen
to be equivalent to the
standard Hille condition \cite{[Hil]}:
\begin{equation}\label{12.15}
\sup_{a>0} \, \, a \, \int_a^\infty V(x) \, dx  < \infty.
\end{equation}
A similar statement is true for the compactness criterion (\ref{12.14}).

\section{Infinitesimal form boundedness}

Maz'ya and Verbitsky \cite{[MV3]} characterized the class of potentials $V \in {\cal D}'(\Bbb{R}^n)$
which are $-\Delta$-form bounded with relative bound zero, i.e.,
for every $\epsilon>0$,
there
exists  $C(\epsilon)>0$ such that
\begin{equation}\label{13.E:1.1}
\left \vert \langle V u, \, u\rangle \right \vert  \leq \epsilon \,
||\nabla u||^2_{L_2(\Bbb{R}^n)} + C(\epsilon) \, ||u||^2_{L_2(\Bbb{R}^n)},
\quad \forall u  \in
C^\infty_0(\Bbb{R}^n).
\end{equation}
In other words, they found necessary and sufficient conditions
for the  {\em infinitesimal form
boundedness\/}  of the potential energy operator $V$ with respect to the
 kinetic energy operator
$H_0 = -\Delta$ on $L_2(\Bbb{R}^n)$.  Here $V$ is an  arbitrary
real- or
complex-valued  potential (possibly a distribution). This notion appeared in relation to the  KLMN theorem
 and has become  an indispensable tool in mathematical
quantum mechanics  and PDE  theory.

The preceding inequality ensures that,
 in case $V$ is real-valued,
 a semi-bounded self-adjoint Schr\"odinger operator $H_V = H_0 + V$
can be defined on $L_2(\Bbb{R}^n)$ so that
the  domain of $Q[u,u]$ coincides
with $W^{1}_2 (\Bbb{R}^n)$. For complex-valued $V$, it follows that
$H_V$ is  an $m$-sectorial operator on $L_2(\Bbb{R}^n)$ with
${\rm Dom}(H_V)\subset W^{1}_2(\Bbb{R}^n)$ (\cite{[RS2]}, Sec. X.2;
\cite{[EE]}, Sec. IV.4).

The characterization of (\ref{13.E:1.1}) found in \cite{[MV3]}
uses only the functions $|\nabla (1-\Delta)^{-1} \, V|$ and
$|(1-\Delta)^{-1} \, V|$, and is based on the
representation:
\begin{equation}\label{13.E:1.rep}
V = {\rm div} \, \vec \Gamma + \gamma, \qquad
\vec \Gamma (x) = -\nabla (1-\Delta)^{-1} \, V, \quad \gamma = (1-\Delta)^{-1} \, V.
\end{equation}
In particular, it is shown that,
necessarily,  $\vec \Gamma \in  L_2(\Bbb{R}^n,loc)^n,
\quad \gamma \in L_1(\Bbb{R}^n,loc),$
and, when $n\ge 3$,
\begin{equation}\label{13.E:1.1b}
\lim_{\delta\to +0} \,
\sup_{x_0\in \Bbb{R}^n} \, \delta^{2-n} \int_{B_\delta(x_0)} \left (
|\vec \Gamma(x)|^2 + |\gamma(x)| \right) \, dx =0,
\end{equation}
once (\ref{13.E:1.1}) holds.

In the opposite direction, it follows from the results in \cite{[MV3]} that (\ref{13.E:1.1})
holds whenever
\begin{equation}\label{13.E:1.1c}
\lim_{\delta\to +0} \,
\sup_{x_0\in \Bbb{R}^n} \, \delta^{2r-n} \int_{B_\delta(x_0)} \left (
|\vec \Gamma(x)|^{2} + |\gamma(x)| \right)^r \, dx =0,
\end{equation}
where $r>1$. Such admissible potentials form  a natural
analog of the Fefferman--Phong class (\ref{e0})  for the infinitesimal
form boundedness problem where cancellations between the positive and negative
parts of $V$  come into play. It
includes functions with
highly oscillatory behavior as well as singular measures, and  properly contains
 the class of potentials  based on the original Fefferman--Phong condition
where  $|V|$ is used  in (\ref{13.E:1.1c}) in place of
$|\vec \Gamma|^{2} + |\gamma|$. Moreover, one can expand
this class further using the sharper condition (\ref{e4})  applied to $|\vec \Gamma|^{2} + |\gamma|$.

A complete  characterization of (\ref{13.E:1.1}) obtained in \cite{[MV3]} is given in the  following theorem
which  provides for deducing
explicit criteria of
 the infinitesimal form
boundedness  in terms of the
{\it nonnegative} locally integrable functions $|\vec \Gamma|^2$ and
$|\gamma|$.

\begin{theorem}\label{13.Theorem II}
 Let $V \in \mathcal {\cal D}'(\Bbb{R}^n)$, $n \ge 2$. The
following statements are equivalent:

{\rm (i)} $V$ is infinitesimally form bounded with respect to $-\Delta$.

{\rm (ii)} $V$ has the form {\rm (\ref{13.E:1.rep})} where
$\vec \Gamma = - \nabla (1-\Delta)^{-1} \, V$,
$\gamma = (1-\Delta)^{-1} \, V$, and the measure $\mu\in M^+(\Bbb{R}^n)$ defined by
\begin{equation}\label{13.E:1.8}
d \mu = \left ( |\vec \Gamma (x)|^2 + |\gamma(x)| \right) \, dx
\end{equation}
has the property that, for every $\epsilon>0$, there exists $C(\epsilon)>0$
such that
\begin{equation}\label{13.E:1.9}
 \int_{\Bbb{R}^n} |u(x)|^2 \, d \mu \le \epsilon \, ||\nabla u||^2_{L_2(\Bbb{R}^n)}
+ C(\epsilon) \, ||\nabla u||^2_{L_2(\Bbb{R}^n)}, \quad \forall u \in C^\infty_0(\Bbb{R}^n).
\end{equation}

{\rm (iii)} For $\mu$ defined by {\rm (\ref{13.E:1.8})},
\begin{equation}\label{13.E:1.11}
 \lim_{\delta \to +0} \, \sup_{P_0: \, {\rm diam} \,  P_0 \le \delta} \,
 \frac {1} {\mu(P_0)} \, \sum_{P \subseteq P_0}  \,  \frac
{\mu(P)^2}  {({\rm mes}_n \, P)^{1- 2/n}}
 = 0,
\end{equation}
where $P$, $P_0$ are dyadic cubes in $\Bbb{R}^n$, i.e., sets of the form
$2^{i} (k + [0, \, 1)^n)$, where $i \in \Bbb{Z}, \, k \in \Bbb{Z}^n$.

{\rm (iv)} For $\mu$ defined by {\rm (\ref{13.E:1.8})},
\begin{equation}\label{13.E:1.10}
\lim_{\delta \to +0} \sup_{F: \, {{\rm diam}} \,  F \le \delta}
\frac {\mu(F)}
{ {\rm{cap}} \, F} = 0,
\end{equation}
where $F$ denotes a  compact set of positive capacity in $\Bbb{R}^n$.

 {\rm (v)} For $\mu$ defined by {\rm (\ref{13.E:1.8})},
\begin{equation}\label{13.E:1.12a}
 \lim_{\delta \to +0} \,  \sup_{x_0\in \Bbb{R}^n} \,  \frac {\left \Vert
 \mu_{B_\delta(x_0)} \, \right \Vert^2_{W^{-1}_2 (\Bbb{R}^n)}}
{ \mu (B_\delta(x_0))}  = 0,
\end{equation}
where  $\mu_{B_\delta(x_0)}$ is the restriction of $\mu$ to the ball
$B_\delta(x_0)$.

{\rm (vi)} For $\mu$ defined by {\rm (\ref{13.E:1.8})},
\begin{equation}\label{13.E:1.13a}
 \lim_{\delta \to +0} \, \sup_{x, \, x_0\in \Bbb{R}^n} \,
 \frac { G_1 \ast \left ( G_1 \ast
\mu_{B_\delta(x_0)} \right )^2 (x)}{G_1 \ast
\mu_{B_\delta(x_0)}(x) } = 0,
\end{equation}
where
 $G_1 \ast \mu = (1-\Delta)^{-\frac 1 2} \mu$ is the Bessel potential of order $1$.
\end{theorem}

 In the one-dimensional case,
 the infinitesimal form boundedness
of the Sturm-Liouville operator $H_V  = -{d^2}/{dx^2} +V$ 
 on
$L_2(\Bbb{R}^1)$
is actually a  consequence of the form boundedness.

\begin{theorem}\label{13.Theorem III}
 Let $V \in \mathcal {\cal D}'(\Bbb{R}^1)$. Then the
following statements
are equivalent.

{\rm (i)} $V$ is infinitesimally form bounded with respect to
$-{d^2}/{dx^2}$.

{\rm (ii)} $V$ is  form bounded with respect to
$-{d^2}/{dx^2}$, i.e.,
$$
|\langle V\, u, \, u\rangle| \le C \, ||u||^2_{W^{1}_2(\Bbb{R}^1)}, \quad
\forall u \in C^\infty_0(\Bbb{R}^1).
$$

{\rm (iii)} $V$ can be represented in the form $V =  {d \Gamma}/ {dx} + \gamma$, where
\begin{equation}\label{13.E:1.12}
\sup_{x \in \Bbb{R}^1} \int_x^{x+1} \left ( |\Gamma(x)|^2 + |\gamma(x)| \right) \,
dx < +\infty.
\end{equation}

{\rm (iv)} Condition {\rm (\ref{13.E:1.12})}
  holds where
$$
\Gamma(x) = \int_{\Bbb{R}^1} {\rm sign} \, (x-t) \, e^{-|x-t|} \, V(t) \, dt, \qquad
\gamma (x) =  \int_{\Bbb{R}^1} e^{-|x-t|} \, V(t) \, dt
$$
are understood in the distributional sense.

{\rm (v)} $V$ belongs to the space $W^{-1}_2(\Bbb{R}^1, unif)$,
with the norm
$$\sup_{x\in \Bbb{R}^1} \|\eta(x-\cdot)\,
V\|_{W^{-1}_2(\Bbb{R}^1)}, \qquad {\rm where} \,\, \eta\in
C^\infty_0(\Bbb{R}^1), \,\, \eta(0) =1.$$

\end{theorem}

The statement   (iii)$\Rightarrow$(i) in Theorem \ref{13.Theorem
III} can be found in \cite{[Sch1]}, Theorem 11.2.1, whereas
(ii)$\Rightarrow$(iv) follows from  results  in \cite{[MV2]},
Theorem 4.2 and \cite{[MV3]},  Theorem 2.5. The equivalence of
(ii) and (v) is proved in \cite{[MSha2]}, Corollary 9.
$\qquad\square$

\section{Kato's condition $K_n$}

Among well-known sufficient conditions for (\ref{13.E:1.1}) 
which ignore possible cancellations, we  mention:
$Q\in L^{\frac n 2} (\Bbb{R}^n) + L^\infty(\Bbb{R}^n)$
($n\ge 3$) and  $Q\in L^{r} (\Bbb{R}^2) + L^\infty(\Bbb{R}^2)$, $r>1$ ($n=2$) (see \cite{[BrK]}),
as well as Kato's condition $K_n$ introduced in \cite{[Ka2]}:
\begin{equation}\label{13.E:K_n}
\lim_{\delta \to +0} \, \sup_{x_0\in\Bbb{R}^n} \,  \int_{B_\delta(x_0)}
 \frac {|Q(x)|} {|x-x_0|^{n-2}}
 \, dx
 = 0, \qquad\qquad\  n \ge 3,
 \end{equation}
 \begin{equation}\label{13.E:K_2}
\lim_{\delta \to +0} \, \sup_{x_0\in\Bbb{R}^n} \,  \int_{B_\delta(x_0)}
\log  \frac 1 {|x-x_0|}  \, |Q(x)|
 \, dx
 = 0, \qquad n =2.
\end{equation}

Kato's class proved to be  especially important in studies of Schr\"odinger semigroups,
Dirichlet forms,
and Harnack  inequalities \cite{[Agm]}, \cite{[ASi]}, \cite{[Sim3]}. Theorem \ref{13.Theorem II} yields
that  (\ref{13.E:1.1}) actually holds for a substantially
broader class of potentials for which $|\vec \Gamma|^2 + |\gamma|
\in K_n$. We emphasize that
no  {\it a priori} assumptions
were imposed  on $C(\epsilon)$ in this theorem.

An observation of Aizenman and Simon
 states that, under the hypothesis
$$C(\epsilon) \le a \, e^{b \, \epsilon^{-p}} \quad  {\rm for} \,\,{\rm some}\,\, a, b>0\,\, {\rm and} \,\,0<p<1,$$
 all   potentials $V$
which obey (\ref{13.E:1.1}) with $|V|$ in place of $V$  are contained in Kato's class.
This was first
proved in \cite{[ASi]}  using the Feynman--Kac formalism. In \cite{[MV3]},
  a sharp
result of this kind is obtained with a simple analytic proof.
 It is shown that
if (\ref{13.E:1.1}) holds with $|V|$ in place of $V \in L^1_{{\rm loc}}(\Bbb{R}^n)$,
then for any $C(\epsilon)>0$,
\begin{equation}\label{13.E:k1}
 \sup_{x_0\in\Bbb{R}^n} \, \int_{B_\delta(x_0)}
 \frac {|V(x)|} {|x-x_0|^{n-2}}
 \, dx \le c \,
\int_{\delta^{-2}}^{+\infty} \frac {\hat C(s)} {s^2} \, ds,
\qquad n \ge 3,
\end{equation}
\begin{equation}\label{13.E:k2}
 \sup_{x_0\in\Bbb{R}^2} \, \int_{B_\delta(x_0)} \log \frac 1 {|x-x_0|} \,
|V(x)|  \, dx \le c \,
\int_{\delta^{-2}}^{+\infty} \frac {\hat C(s)} {s^2 \log s} \, ds,
\qquad n =2,
\end{equation}
where $c$ is a constant which depends only on $n$, and $\delta$ is sufficiently small. Here
$\hat C(s) = \inf_{\epsilon>0} \, \{C(\epsilon) + s \, \epsilon\}$ is the Legendre transform of $-C(\epsilon)$.
In particular, it follows that the condition $C(\epsilon) \le a \, e^{b \, \epsilon^{-p}}$ for  {\it any}
 $p>0$ is enough  to ensure
that $V \in K_2$ in the more subtle two-dimensional case.

\section{Trudinger's subordination for the Schr\"odinger operator}

In \cite{[MV3]} inequality (\ref{13.E:1.1}) is
 studied also under the assumption that  $C(\epsilon)$  has   power growth,
i.e., there exists $\epsilon_0>0$ such that
  \begin{equation}\label{13.E:1.2}
\left \vert \langle V u, \, u\rangle \right \vert  \leq \epsilon \,
||\nabla u||^2_{L_2(\Bbb{R}^n)} + c \,   \epsilon^{-\beta}  \,
||u||^2_{L_2(\Bbb{R}^n)},  \quad \forall u
\in C^\infty_0(\Bbb{R}^n),
\end{equation}
 for every
$\epsilon\in(0, \epsilon_0)$, where $\beta>0$. Such inequalities
appear in  studies of elliptic PDE with measurable
coefficients \cite{[Tru]}, and have been  used extensively in spectral theory
 of the  Schr\"odinger operator  \cite{[ASi]},  \cite{[Ka2]}, \cite{[RS2]}, \cite{[RSS]}, \cite{[Sch1]}, \cite{[Dav1]}, \cite{[Dav2]}, \cite{[LPS]}, \cite{[Sim3]}.

As it turns out, it is still possible to
characterize (\ref{13.E:1.2}) using only $|\vec \Gamma|$ and
$|\gamma|$ defined by (\ref{13.E:1.rep}), provided $\beta>1$. It is shown in \cite{[MV3]} that
in this case (\ref{13.E:1.2})
 holds if and only if
both of the following conditions hold:
\begin{equation}\label{13.E:1.2a}
 \sup_{{x_0\in \Bbb{R}^n}\atop{0<\delta< \delta_0}} \,
\delta^{2 \frac {\beta-1}{\beta+1}-n}
 \int_{B_\delta(x_0)}
|\vec \Gamma(x)|^{2} \, dx < +\infty,
\end{equation}
\begin{equation}\label{13.E:1.2b}
\sup_{{x_0\in \Bbb{R}^n}\atop{0<\delta< \delta_0}} \,
\delta^{\frac {2\beta}{\beta+1}-n}
 \int_{B_\delta(x_0)}  |\gamma(x)| \, dx < +\infty,
\end{equation}
for some $\delta_0>0$.
However, in the case $\beta\le 1$
this is no longer true. For $\beta=1$, (\ref{13.E:1.2a}) has to be replaced with
the condition that $\vec \Gamma$ is in the local ${\rm BMO}$ space,
or respectively is H\"older-continuous of order ${(1-\beta)}/{(1+\beta)}$ if
$0<\beta<1$.

In the homogeneous case $\epsilon_0=+\infty$, (\ref{13.E:1.2})
is equivalent to the  {\em multiplicative inequality\/} (\ref{x1})
where $\tau = {\beta}/{(1 + \beta)} \in (0, \, 1)$. In
spectral theory,  (\ref{x1})
 is  referred to as the form $\tau$-subordination property
(see, e.g., \cite{[Agr]}, \cite{[Gr1]}, \cite{[Gr2]}, \cite{[MM1]},  \cite{[RSS]}, Sec. 20.4).

For {\em nonnegative\/}  potentials $V$,  where $V$ coincides with a
locally finite measure $\mu$  on $\Bbb{R}^n$,  inequality
(\ref{x1}) is equivalent to (\ref{x2}) \cite{[Maz7]}, Sect. 1.4.7.
 For general $V$,  the following result is obtained in \cite{[MV3]}. If $\tau> 1/ 2$,
then  {\rm (\ref{x1})}
holds if and only if $\nabla \Delta^{-1} V$ lies in  the
 Morrey  space
$\mathcal{L}^{2, \, \lambda}(\Bbb{R}^n)$, where $\lambda= n +2 -4\tau$. The  Morrey space $\mathcal{L}^{r, \, \lambda}(\Bbb{R}^n)$ ($r>0,  \lambda>0$)
consists of $f \in  L^r (\Bbb{R}^n, loc)$ such that
$$
\sup_{x_0\in \Bbb{R}^n, \, \delta>0} \,\,  \frac
{1}{|B_\delta(x_0)|^{\lambda/n}} \int_{B_\delta(x_0)}
|f(x)|^r \, dx < + \infty.
$$
We refer to \cite{[Pe]} for the overview of Morrey spaces.
For $\tau=1/ 2$,  (\ref{x1}) holds if and only if
$\nabla \Delta^{-1} V \in {\rm BMO}(\Bbb{R}^n)$, and for
$0<\tau<1/ 2$, whenever  $\nabla \Delta^{-1} V$ is in the H\"older  class $C^{1-2\tau}(\Bbb{R}^n)$.
These different
characterizations  are equivalent
to  {\rm (\ref{x2})}
 if $V$ is a nonnegative measure.

\section{Discreteness of the spectrum of $-\Delta +\Bbb{V}$ with nonnegative potential}

Friedrichs  proved in 1934 \cite{[Fr]} that to the left of the point $\mathop{\hbox{{\rm lim inf}}}_{|x|\to\infty} \, V(x)$,  
the spectrum of the Schr\"odinger operator $H_V=-\Delta+V$ in
$L_2(\Bbb{R}^n)$ with a locally integrable potential $V$ is
discrete,  and hence all spectrum is discrete provided
$V(x)\to+\infty$ as $|x|\to\infty$. By Rellich's criterion
\cite{[Re2]}, the spectrum of $H_V$ is purely discrete if and only
if the ball $(H_V u, u)\leq 1$ is a compact subset of
$L_2(\Bbb{R}^n)$. This and
   the discreteness of spectrum easily imply that for every $d>0$
\begin{equation}\label{h1}
\int_{Q_d}V(x)dx\to +\infty \quad {\rm as} \quad Q_d\to\infty,
\end{equation}
where $Q_d$ is a closed cube with the edge length $d$ and with the
edges parallel to coordinate axes,
$Q_d\to\infty$ means that the cube $Q_d$ goes to infinity (with fixed $d$). In 1953 
A.M.Molchanov \cite{[Mol]} proved that this condition is in fact
necessary and sufficient in case $n=1$ but not sufficient for $n\ge 2$.
He also discovered a modification of (\ref{h1}) which is
equivalent to the discreteness of spectrum in the case  $n\ge 2$:
\begin{equation}\label{h2}
\inf_F\int_{Q_d\setminus F} V(x) dx \to +\infty \quad {\rm as} \quad
Q_d\to\infty,
\end{equation}
where infimum is taken over all compact subsets $F$ of $Q_d$ which are called
{\it negligible}.  The negligibility of $F$ in the sense of Molchanov
means that
${\rm cap} F\le\gamma\, {\rm cap}(Q_d)$, where ${\rm cap}$ is the Wiener capacity and $\gamma>0$
is a sufficiently small constant. More precisely, Molchanov proved
that one can take $\gamma=c_n$
where $\gamma=c_n=(4n)^{-4n}({\rm cap}(Q_1))^{-1}$ for $n\ge 3$.

As early as in 1953, I.M.Gelfand raised the question about
the best possible constant $c_n$ (personal communication). We describe  results from \cite{[MShu1]}, where a complete answer to this question is given.

Let $\Bbb{V}$ be a positive Radon measure in  an open set $\Omega\subset\Bbb{R}^n$, absolutely continuous with respect to the Wiener capacity. We will consider
the Schr\"odinger operator $H_{\Bbb{V}}$ which is formally given by the expression $-\Delta+\Bbb{V}$. It
is defined in $L_2(\Omega)$ by the closure of the quadratic form $Q[u,u]$ with the domain
  $C_0^\infty(\Omega)$.

Instead of the cubes $Q_d$, a more
general family of test bodies will be used. Let us start with a standard open set ${\cal G}\subset \Bbb{R}^n$.
We assume that ${\cal G}$ satisfies the following conditions:

(a) ${\cal G}$ is bounded and star-shaped  with respect
to the  ball $B_\rho$;

(b) ${\rm diam}({\cal G})=1$.

For any positive $d>0$ denote by ${\cal G}_d(0)$ the body
$\{x: \;d^{-1}x\in{\cal G}\}$
which is  homothetic to ${\cal G}$ with coefficient $d$  and with the center of homothety at $0$.
We will denote by ${\cal G}_d$ a body which is obtained from ${\cal G}_d(0)$ by a parallel
translation: ${\cal G}_d(y)=y+{\cal G}_d(0)$ where
$y$ is an arbitrary  vector in $\Bbb{R}^n$.

The notation ${\cal G}_d\to\infty$ means that
the distance from ${\cal G}_d$ to $0$ goes to infinity.

\noindent
{\bf Definition.}
Let $\gamma\in(0,1)$. The {\it negligibility class} ${\cal N}_{\gamma}({\cal G}_d;\Omega)$ consists
of all compact sets $F\subset\overline{\cal G}_d$ satisfying the following conditions:
\begin{equation}\label{E:F-inclusion}
\bar{\cal G}_d\setminus\Omega\subset F\subset \bar {\cal G}_d\;,
\end{equation}
and
\begin{equation}\label{E:G-negligibility}
{\rm cap}\, F\le \gamma\, {\rm cap}\, \bar{G}_d.
\end{equation}

Now we formulate the main result of the work in \cite{[MShu1]}
about the discreteness of spectrum.

\begin{theorem}\label{T:discr}

{\rm (i)} {\rm (Necessity)} Let the spectrum of $H_\Bbb{V}$ be discrete. Then for
every  function $\gamma:(0,+\infty)\to (0,1)$ and every $d>0$
\begin{equation}\label{E:inf-cond}
\inf_{F\in{\cal N}_{\gamma(d)}({\cal G}_d,\Omega)}\; \Bbb{V}({\bar {\cal G}_{d}\setminus F)}\to +\infty
\quad {\rm as} \quad {\cal G}_d\to\infty.
\end{equation}

{\rm (ii)} {\rm(Sufficiency)}
Let a function $d\to\gamma(d)\in(0,1)$ be defined for $d>0$
in a neighborhood of $0$, and satisfy
\begin{equation}\label{h5}
\mathop{\hbox{{\rm lim sup}}}_{d\downarrow 0} d^{-2}\, \gamma(d) = +\infty.
\end{equation}
Assume that
there exists $d_0>0$
such that(\ref{E:inf-cond}) holds for every $d\in (0,d_0)$.
Then the spectrum of $H_\Bbb{V}$ in $L_2(\Omega)$ is discrete.

\end{theorem}

Let us make some comments about this theorem.

\begin{remark}\label{R:sequence}
{\rm
 It suffices for the discreteness of spectrum of $H_{\Bbb{V}}$
that the condition (\ref{E:inf-cond}) holds
only for a sequence of $d$'s, i.e. $d\in\{d_1,d_2,\dots\}$, $d_k\to 0$ and
$d_k^{-2}\gamma(d_k)\to +\infty$ as $k\to +\infty$.
}
\end{remark}

\begin{remark}\label{R:other-suff}
{\rm
The condition (\ref{E:inf-cond}) in the sufficiency part can be
replaced  by a weaker requirement: there exist $c>0$ and $d_0>0$
such that for every $d\in (0,d_0)$
there exists $R>0$ such that
\begin{equation}\label{E:inf-ineq}
d^{-n}\inf_{F\in{\cal N}_{\gamma(d)}({\cal G}_d,\Omega)}\; \Bbb{V}({\bar{\cal G}_{d}\setminus F)}
 \ge cd^{-2}\gamma(d),
\end{equation}
whenever $\bar{\cal G}_d\cap(\Omega\setminus B_R)\ne\emptyset$
(i.e. for distant bodies ${\cal G}_d$ having non-empty intersection with $\Omega$).
Moreover, it suffices that the condition
(\ref{E:inf-ineq}) is satisfied for a sequence $d=d_k$ satisfying the condition
formulated in Remark \ref{R:sequence}.

Note that unlike (\ref{E:inf-cond}), the condition (\ref{E:inf-ineq}) does not require
that the left-hand side goes to $+\infty$ as ${\cal G}_d\to\infty$. What is actually needed
is that the left-hand side has a certain lower bound, depending on $d$ for arbitrarily small
$d>0$ and distant test bodies ${\cal G}_d$. Nevertheless,
the conditions (\ref{E:inf-cond}) and (\ref{E:inf-ineq})
are equivalent because each  is equivalent to the discreteness of spectrum.
}
\end{remark}

\begin{remark}\label{R:Molchanov}
{\rm
If we take $\gamma=const\in(0,1)$, then Theorem \ref{T:discr} gives Molchanov's result, but
with the constant $\gamma=c_n$
replaced by an arbitrary constant $\gamma\in (0,1)$. So Theorem \ref{T:discr} contains an answer
to the above-mentioned Gelfand's question.
}
\end{remark}

\begin{remark}\label{R:ga-equiv}
{\rm
For any two functions $\gamma_1,\gamma_2:(0,+\infty)\to (0,1)$ satisfying the requirement
(\ref{h5}), the conditions (\ref{E:inf-cond}) are equivalent, and so are
the conditions (\ref{E:inf-ineq}), because any of these conditions is equivalent
to the discreteness of spectrum.

It follows that the conditions  (\ref{E:inf-cond}) for different constants $\gamma\in(0,1)$
are equivalent. In the particular case, when the measure $\Bbb{V}$ is absolutely continuous
with respect to the Lebesgue measure, we see that the conditions (\ref{h2})
with different constants $\gamma\in(0,1)$ are equivalent.
}
\end{remark}

\begin{remark}\label{R:domains}
{\rm The results above are new even for the operator $-\Delta$
in $L_2(\Omega)$ (for an arbitrary open set $\Omega\subset
\Bbb{R}^n$ with the Dirichlet boundary conditions on
$\partial\Omega$). In this case the discreteness of spectrum is
completely determined  by the geometry of $\Omega$. More
precisely, for the discreteness of spectrum of $H_0$ in
$L_2(\Omega)$ it is necessary and sufficient that there exists
$d_0>0$ such that for every $d\in (0,d_0)$
\begin{equation}\label{E:omega-cond}
\mathop{\hbox{{\rm lim inf}}}_{{\cal G}_d\to\infty}{\rm cap}(\bar{\cal G}_d\setminus \Omega)\ge \gamma(d)\, {\rm cap}\, \bar{\cal G}_d,
\end{equation}
where $d\to \gamma(d)\in (0,1)$ is a function, which is  defined in a neighborhood
of $0$ and satisfies (\ref{h5}). The conditions (\ref{E:omega-cond})
with different functions $\gamma$, satisfying the conditions above, are equivalent.
This is a non-trivial property of capacity.
It is necessary for the discreteness of spectrum  that (\ref{E:omega-cond})
holds for every function $\gamma:(0,+\infty)\to (0,1)$ and every $d>0$, but this condition
may not be sufficient if $\gamma$ does not satisfy (\ref{h5})
(see Theorem \ref{T:precise} below).
}
\end{remark}

The following result demonstrates that the condition (\ref{h5}) is precise.

\begin{theorem}\label{T:precise}
Assume that $\gamma(d)=O(d^2)$ as $d\to 0$. Then there exists an open set
$\Omega\subset\Bbb{R}^n$ and $d_0>0$ such that for every $d\in(0,d_0)$
the condition  (\ref{E:omega-cond}) is satisfied
but the spectrum of $-\Delta$ in $L_2(\Omega)$ with the Dirichlet boundary conditions
is not discrete.
\end{theorem}

\section{Strict positivity of the spectrum of $-\Delta + \Bbb{V}$}

We  say that the operator $H_{\Bbb{V}}$, the same as in Sect. 11,
is {\it strictly positive} if its spectrum does not contain $0$.
Equivalently, we can say that the spectrum is separated from $0$.
 The strict positivity is equivalent
to the existence of $\lambda>0$ such that
\begin{equation}\label{E:form-pos}
Q[u,u]\ge \lambda \|u\|^2_{L_2(\Omega)},\quad u\in C_0^\infty(\Omega).
\end{equation}

The characterization of positivity of the spectrum in the next
theorem with Molcha- nov's negligible sets in the formulation was
found in \cite{[Maz5]} (see also \cite{[Maz7]}, Sect. 12.5). The
present stronger version is obtained by Maz'ya and Shubin in
\cite{[MShu1]}.

\begin{theorem}\label{T:positivity}
{\rm (i) (Necessity)} Let us assume that $H_{\Bbb{V}}$ is strictly positive, so that
(\ref{E:form-pos}) is satisfied with a constant $\lambda>0$. Let us take an arbitrary
$\gamma\in (0,1)$.
Then there exist $d_0>0$ and $\varkappa>0$ such that
\begin{equation}\label{E:pos-cond2}
d^{-n}\inf_{F\in {\cal N}_{\gamma}({\cal G}_d,\Omega)} \Bbb{V}(\bar{\cal G}_d\setminus F)\ge \varkappa
\end{equation}
for every $d>d_0$ and every ${\cal G}_d$.

{\rm (ii) (Sufficiency)} Assume that there exist $d>0$, $\varkappa>0$ and $\gamma\in (0,1)$,
such that (\ref{E:pos-cond2}) is satisfied for every ${\cal G}_d$.
Then the operator $H_{\Bbb{V}}$ is strictly positive.

Instead of all bodies ${\cal G}_d$ it is sufficient to take only the ones from
a finite multiplicity covering (or tiling) of $\Bbb{R}^n$.
\end{theorem}

\begin{remark}\label{R}
{\rm
Considering the Dirichlet Laplacian  in $L_2(\Omega)$ we see from Theorem \ref{T:positivity} that for any choice of a constant $\gamma\in (0, 1)$ and a standard body ${\cal G}$, the strict positivity of $-\Delta$ is equivalent to the following condition
\begin{equation}\label{k9x}
\exists\, d>0, \,\, {\rm such}\,\, {\rm that}\,\, {\rm cap} (\bar{\cal G}_d\cap (\Bbb{R}^n\backslash\Omega)) \geq \gamma\, {\rm cap}(\bar{\cal G}_d) \,\, {\rm for}\,\, {\rm all}\,\, {\cal G}_d.
\end{equation}

\noindent
In particular, it follows that for two different $\gamma$'s  these conditions are equivalent. Noting that $\Bbb{R}^n \backslash\Omega$ can be an arbitrary closed subset in $\Bbb{R}^n$, we get a property of the Wiener capacity, which is obtained as a byproduct of our spectral theory arguments.
}
\end{remark}

\section{Two-sided estimates for the bottom of spectrum and essential spectrum}

Let $\lambda= \lambda(\Omega, H_{\Bbb{V}})$ denote the greatest lower bound of the spectrum of the Schr\"odinger operator $H_{\Bbb{V}}$ handled in the two preceding sections. By ${\cal Q}_\Omega$ we denote the set of all cubes $Q_d$ having a negligible intersection with the complement of $\Omega$ (in Molchanov's sense), i.e.
\begin{equation}\label{j1}
{\cal Q}_\Omega=\{Q_d: \; {\rm cap}(Q_d\setminus \Omega)\le\gamma\,{\rm cap}\, Q_d\},
\end{equation}
where $\gamma$ is sufficiently small.

For $\Bbb{V}=0$ Maz'ya \cite{[Maz6]} and for the general case 
Maz'ya and Otelbaev \cite{[MaO]} (see also \cite{[Maz7]}, Ch.\,12) 
 obtained the two-sided estimate
\begin{equation}\label{1.6}
c_1D^{-2}\le \lambda \le c_2D^{-2},
\end{equation}
with $D=D(\Omega,\Bbb{V})$  given by
\begin{equation}\label{1.7}
D=\sup_{Q_d\in{\cal Q}_\Omega}\left\{d:\;d^{n-2}\ge \inf_F\Bbb{V}(Q_d\backslash F)
\right\}.
\end{equation}
Here  $c_1$  and $c_2$ are positive  constants which depend only upon $n$  and
the infimum is taken over all  sets $F$,
satisfying
$$Q_d\setminus \Omega\subset F\subset Q_d \quad {\rm and}\quad {\rm cap}F\le\gamma\,{\rm cap}\, Q_d.
$$
Clearly, the strict
positivity
of $H_{\Bbb{V}}$ is equivalent to the condition $D<\infty$. 

In particular,
if $\Bbb{V}= 0$, then $D$ becomes the capacitary interior diameter of $\Omega$:
the maximal size $d$ of cubes $Q_d$ with the negligible intersection with
$\Bbb{R}^n\setminus\Omega$ (i.e. cubes $Q_d\in{\cal Q}_\Omega$).  The notion of the capacitary interior diameter was introduced in \cite{[Maz6]}.
Explicit estimates for $c_1$ and $c_2$ in the case $\Bbb{V}=0$ are given in \cite{[MShu2]}, where also the smallness condition of $\gamma$ is replaced with $\gamma\in (0,1)$.

The paper \cite{[MaO]} (see also \cite{[Maz7]}, Ch. 12)  also contains a two-sided estimate
for the bottom of the essential spectrum of $H_{\Bbb{V}}$ in $L_2(\Omega)$. We will denote
this bottom  by $\Lambda=\Lambda(\Omega;H_{\Bbb{V}})$. The estimate has the form
\begin{equation}\label{1.8}
c_1D_\infty^{-2}\le \Lambda\le c_2 D_\infty^{-2},
\end{equation}
where
\begin{equation}\label{1.9}
D_\infty=\lim_{R\to\infty} D(\Omega\setminus \bar B_R,\Bbb{V},\gamma).
\end{equation}

Note that the discreteness of spectrum of $H_{\Bbb{V}}$ is equivalent to the equality
$\Lambda=+\infty$. Therefore, (\ref{1.8}) implies that  $D_{\infty}=0$
is  necessary and sufficient  for the pure discreteness of  the spectrum of $H_{\Bbb{V}}$.

In the recent work \cite{[Tay]} M. Taylor found another criterion of  discreteness of the spectrum of $H_{\Bbb{V}}$ stated in terms of the so-called scattering length of $\Bbb{V}$.

\section{Structure of the essential spectrum of $H_{\Bbb{V}}$}

 We use the
same notation as in Sections 16--18. The following result is due to Glazman \cite{[Gl]}.

\begin{lemma}\label{Lem3}
If the spectrum of $H_{\Bbb{V}}$ is not purely discrete, the essential spectrum of $H_{\Bbb{V}}$ extends to infinity. Moreover, if $0$ belongs to the essential spectrum of $H_{\Bbb{V}}$, then this spectrum coincides with $[0,\infty)$.
\end{lemma}

{\bf Proof.} Let $\Lambda$ be the bottom of the essential spectrum. Then there exists a sequence of real-valued functions $\{\varphi_\nu\}_{\nu\geq 1}$ in $C^\infty_0(\Omega)$ subject to the conditions
\begin{equation}\label{l1}
\|\varphi_\nu\|_{L_2(\Omega)} = 1, \qquad \varphi_\nu \to 0\quad {\rm weakly}\,\, {\rm in} \,\, L_2(\Omega),
\end{equation}
\begin{equation}\label{l2}
\|(H_{\Bbb{V}} - \Lambda)\, \varphi_\nu\|_{L_2(\Omega)} \to 0.
\end{equation}

We set
$$u_\nu =\varphi_\nu{\rm exp}\Bigl( i(\alpha - \Lambda)^{1/2}\sum_{k=1}^n x_k\Bigr),$$
where $\alpha >\Lambda$. We see that $u_\nu$ satisfies (\ref{l1}) and that
$$
\|(H_{\Bbb{V}} - \alpha)\, u_\nu\|_{L_2(\Omega)}^2 = \|(H_{\Bbb{V}} - \Lambda)\, \varphi_\nu\|_{L_2(\Omega)}^2 + 4(\alpha - \Lambda) \|\sum _{k=1}^n \partial \varphi_\nu /\partial x_k\|_{L_2(\Omega)}^2.$$
Since the right-hand side does not exceed
$$\|(H_{\Bbb{V}} - \Lambda)\, \varphi_\nu\|_{L_2(\Omega)}^2 + 4(\alpha - \Lambda) \, Q[\varphi_\nu, \varphi_\nu],$$
we have
$$\|(H_{\Bbb{V}} - \alpha)\, u_\nu\|_{L_2(\Omega)}^2 \leq
\|(H_{\Bbb{V}} - \Lambda)\, \varphi_\nu\|_{L_2(\Omega)}^2$$
$$ + 4n(\alpha - \Lambda)\, \|(H_{\Bbb{V}} - \Lambda)\, \varphi_\nu\|_{L_2(\Omega)} + 4n\Lambda(\alpha - \Lambda).$$
By (\ref{l2})
$$\mathop{\hbox{{\rm lim sup}}}_{\nu\to \infty} \|(H_{\Bbb{V}} - \alpha)\, u_\nu\|_{L_2(\Omega)}\leq \rho(\alpha),$$
where $\rho(\alpha) = 2(n\, \Lambda(\alpha - \Lambda))^{1/2}$. It
follows that  any segment $[\alpha - \rho(\alpha), \alpha
+\rho(\alpha)]$ contains points of the essential spectrum. If, in
particular, $\Lambda = 0$ then every positive $\alpha$ belongs to
the essential spectrum. $\qquad\square$

In concert with this lemma the pairs $(\Omega, \Bbb{V})$ can be divided into three non-over- lapping  classes.
The first class includes $(\Omega, \Bbb{V})$ such that the spectrum of $H_{\Bbb{V}}$ is discrete.
The pair $(\Omega, \Bbb{V})$ belongs to the second class if the essential spectrum of $H_{\Bbb{V}}$
is unbounded and strictly positive. Finally, $(\Omega, \Bbb{V})$ is of the third class if
the essential spectrum of $H_{\Bbb{V}}$ coincides with $[0,\infty)$.

By (\ref{1.7}) and (\ref{1.8}) the three classes can be described as follows:

(i) $(\Omega, \Bbb{V})$ belongs to the first class if and only if $D_\infty =0$.

(ii) $(\Omega, \Bbb{V})$ belongs to the second class if and only if $D_\infty >0$.

(iii) $(\Omega, \Bbb{V})$ belongs to the third class if and only if $D =\infty$ (or, equivalently, $D_\infty =\infty$).

By Theorem \ref{T:discr}, the condition (i) is equivalent to
(\ref{E:inf-cond}) and (ii) holds if and only if
(\ref{E:pos-cond2}) is valid. Finally,  (iii) is equivalent to the
failure of (\ref{E:omega-cond}) by Theorem \ref{T:positivity} and
Lemma \ref{Lem3}. In other words, (iii) holds if and only if
$$\mathop{\hbox{{\rm lim inf}}}_{{\cal G}_d \to \infty}
\inf_{F\in{\cal N}_{\gamma(d)}({\cal G}_d,\Omega)}\; \Bbb{V}({\bar
{\cal G}_{d}\setminus F)} =0$$
for every $d>0$.

There are still open problems in the study of conditions insuring the
absolute continuity of the spectrum. B.\,Simon conjectured in \cite{[Sim5]}  that the absolutely continuous spectrum of $H_V$ with $V$ subject to
$$\int_{\Bbb{R}^n} (V(x))^2{(1+ |x|)^{1-n}}\, dx<\infty$$
is the positive real axis. An interesting sufficient condition supporting this conjecture:  
$$V = {\rm  div} \, \vec \Gamma \quad  {\rm and} \quad |V(x)| + |\vec \Gamma(x)| \leq c\, (1+ |x|)^{-1/2-\varepsilon}, \,\,\, \varepsilon>0,$$
has been obtained recently by Denisov in \cite{[Den]}.

\section{Two measure boundedness and compactness criteria}

The following lemma is a particular case of a more general result
from \cite{[Maz4]} (see also \cite{[Maz7]}, Sect. 2.3.7).

\begin{lemma}\label{Lem99}
Let $\mu$ and $\nu$ be two nonnegative Radon measures in $\Omega\subset \Bbb{R}^n$, $n\geq 1$.
The inequality
\begin{equation}\label{t99}
\int_\Omega |u|^2\, d\mu \leq C\Bigl( \int_\Omega |\nabla u|^2\, dx + \int_\Omega |u|^2\, d\nu\Bigr)
\end{equation}
holds for all $u\in C^\infty_0(\Omega)$ if and only if there
exists a constant $K>0$ such that for all open bounded sets $g$
and $G$ subject to ${\bar g}\subset G$, ${\bar G}\subset \Omega$,
the inequality
\begin{equation}\label{s99}
\mu(g) \leq K({\rm cap}_G\, {\bar g} + \nu(G) )
\end{equation}
holds.
\end{lemma}
This criterion becomes more transparent in the one-dimensional case because the capacity admits an explicit
representation.

\begin{theorem}\label{Th.t99}
{\rm (see \cite{[Maz9]})}
Let $n=1$ and let $\sigma_d$ denote the open interval $(x-d, x+d)$. Inequality
\begin{equation}\label{v99}
\int_\Omega |u|^2\, d\mu \leq C\Bigl( \int_\Omega |u'|^2\, dx + \int_\Omega |u|^2\, d\nu\Bigr)
\end{equation}
holds for all $u\in C^\infty_0(\Omega)$ if and only if
\begin{equation}\label{y99a}
\mu(\sigma_d(x)) \leq const (\tau^{-1} + \nu(\sigma_{d+\tau}(x))),
\end{equation}
where $x$, $d$, and $\tau$ are such that ${\overline {\sigma_{d+\tau}(x)}}\subset\Omega$, is valid without complementary assumptions about $\mu$ and $\nu$. The sharp constant $C$ in (\ref{v99}) is equivalent to
$$\sup_{x,d,\tau}\frac{\mu(\sigma_d(x))}{\tau^{-1} + \nu(\sigma_{d+\tau}(x))},$$
where $x$, $d$, and $\tau$ are the same as in (\ref{y99a}).
\end{theorem}

We define the space $\ring W^1_2(\nu)$ as the closure of $C^\infty_0(\Omega)$ with respect to the norm
$$\|f\|_{\ring W^1_2(\nu)} = \Bigl(\int_\Omega |f'(x)|^2\, dx + \int_\Omega |f(x)|^2\, d\nu\Bigr)^{1/2}.$$
The condition (\ref{y99a}) is a criterion of boundedness for the
embedding operator $I: \ring W^1_2(\nu) \to L_2(\mu)$.

The next theorem contains
 a two-sided estimate for the essential norm of $I$. We recall that the essential norm of a bounded
 linear operator $A$ acting from $X$ into $Y$, where $X$ and $Y$ are linear normed spaces, is defined by
 $${\rm ess}\, \|A\| = \inf_T \|A-T\|$$
 with infimum taken over all compact operators $T: X\to Y$.

 \begin{theorem}\label{Th.v99}
 {\rm \cite{[Maz9]}} Let
 $$E(\mu,\nu): \lim_{M\to \infty} \sup_{x,d,\tau}\frac{\mu(\sigma_d(x)\backslash [-M,M])}
 {\tau^{-1} + \nu(\sigma_{d+\tau}(x))}.$$
 There exist positive constants $c_1$ and $c_2$ such that
 \begin{equation}\label{y99}
 c_1\, E(\mu,\nu)^{1/2} \leq {\rm ess} \|I\|\leq c_2\, E(\mu,\nu)^{1/2}.
 \end{equation}

 In particular, the operator $I$ is compact if and only if
 $$\lim_{M\to\infty}\sup_{x,d,\tau}\frac{\mu(\sigma_d(x)\backslash [-M,M])}
 {\tau^{-1} + \nu(\sigma_{d+\tau}(x))} =0,$$
 where $x$, $d$, and $\tau$ are the same as in (\ref{y99a}).
\end{theorem}

\end{document}